\documentclass[a4paper,12pt]{article}
\usepackage[dvips]{graphicx}

\usepackage{pst-grad} 

\usepackage{tikz-cd}
\usetikzlibrary{matrix,arrows,decorations.pathmorphing}

\usepackage{cmdtrack}
\usepackage{amsthm}
\usepackage{url}

\def\Nt{\tilde{\cal N}}
\def\N{{\cal N}}
\def\A{{\cal A}}

\def\G{{\cal G}}
\def\St{\tilde{\Sigma}}
\def\Ss{\tilde{\cal S}}
\def\Sb{\breve{\cal S}}
\def\Et{\tilde{\cal E}}
\def\C{{\cal C}}
\def\F{{\cal F}}
\def\H{{\cal H}}

\def\Bl{\mbox{\rm Bl}}
\def\conv{\mbox{\rm conv}}
\def\st{\mbox{\rm st}}
\def\int{\mbox{\rm int}}
\def\aff{\mbox{\rm aff}}
\def\tr{\mbox{\rm tr}}

\def\cirk{\,{\raisebox{.3ex}{\tiny $\circ$}}\,}

\def\pl{\!+\!}
\def\mn{\!-\!}

\def\prop#1#2{\vspace{2ex} \noindent{\sc #1.} {\it #2} \par \vspace{2ex}}
\def\propo#1#2{\vspace{2ex} \noindent{\sc #1.} {\it #2}}
\def\oprop#1#2{\noindent{\sc #1.} {\it #2} \par \vspace{2ex}}
\def\dkz{\noindent{\sc Proof. }}
\def\qed{\hfill $\dashv$\vspace{2ex}}

\def\ccup{\raisebox{-.5pt}{$\cup$}}
\def\R{{\mathbf R}}

\begin{document}

\title{\textbf{On stretching the interval simplex-permutohedron}}
\author{{\sc Zoran Petri\' c}
\\[1ex]
{\small Mathematical Institute, SANU}\\[-.5ex]
{\small Knez Mihailova 36, p.f.\ 367, 11001 Belgrade,
Serbia}\\[-.5ex]
{\small email: zpetric@mi.sanu.ac.rs}}
\date{}
\maketitle

\begin{abstract}
\noindent A family of polytopes introduced by E.M.\ Feichtner, A.\
Postnikov and B.\ Sturmfels, which were named nestohedra, consists
in each dimension of an interval of polytopes starting with a
simplex and ending with a permutohedron. This paper investigates a
problem of changing and extending the boundaries of these
intervals. An iterative application of Feichtner-Kozlov procedure
of forming complexes of nested sets is a solution of this problem.
By using a simple algebraic presentation of members of nested sets
it is possible to avoid the problem of increasing the complexity
of the structure of nested curly braces in elements of the
produced simplicial complexes.
\end{abstract}

\vspace{.5ex}


\section{Introduction}
A history of the family of polytopes that appear in the title of
this paper under the name ``the interval simplex-permutohedron'',
may be traced to Appendix~B of \cite{S97b}, where Stasheff and
Shnider gave a procedure of truncation of an $(n\mn
2)$-dimensional simplex, which realizes the associahedron $K_n$ as
a convex polytope. This procedure implicitly uses the fact that
$K_n$, as an abstract polytope, is induced by (constructions of)
the graph whose edges connect the neighbour symbols of binary
operations in a term with $n$ variables. It is also noted in
\cite{S97b} that this procedure could be modified in order to
obtain cyclohedra out of simplices; the edges of graphs that
induce cyclohedra would connect the neighbour symbols of binary
operations in cyclic terms.

Carr and Devadoss used graphs in \cite{CD06} more explicitly to
define a family of polytopes that includes simplices,
associahedra, cyclohedra, permutohedra and many others. This
approach, which is based on the concept of tubings is further
developed in \cite{DF08}, \cite{D09} and \cite{FS09}.

The first definition of the notion of nested sets, which plays an
important role in this paper, is given by Fulton and MacPherson in
\cite{FMP94}. De Concini and Procesi used this notion and
introduced the notion of building sets in \cite{DeCP95}. Non
linear nestohedra are obtained already in \cite{Ga03} and
\cite{Ga04} by Gaiffi in the theory of real De Concini-Procesi
models by gluing several manifolds with corners.

Feichtner and Kozlov defined in \cite{FK04} the notions of
building sets and nested sets in a pure combinatorial manner in
order to give an abstract framework for the incidence
combinatorics occurring in some constructions in algebraic
geometry (see for example \cite{DeCP95}). Feichtner and Sturmfels
in \cite{FS05} and Postnikov in \cite{P09} used building sets
independently, in a more narrow context, to describe the face
lattices of a family of simple polytopes named nestohedra in
\cite{PRW08}.

An alternative, inductive, approach that leads to the same family
of polytopes is given in \cite{DP10c}. This was an independent
discovery motivated by some related work done in \cite{DP06},
\cite{DP10a} and \cite{DP10b}.

For a fixed dimension $n$, the family of polytopes in question may
be considered as an interval starting with an $n$-dimensional
simplex and ending with an $n$-dimensional permutohedron, while
the partial ordering of the interval is induced by the relation
$\subseteq$ on the corresponding building sets. Such intervals
contain some polytopes, like for example associahedra and
cyclohedra, which are widely used in topology and category theory.
This is the reason why we give the advantage to the families of
simple polytopes (though the corresponding families of simplicial
polytopes are easier to handle). However, if one prefers to
simplify everything a bit, the mathematical content of the results
would not change if we switch to the intervals starting with an
$n$-dimensional simplex and ending with the polytope polar to an
$n$-dimensional permutohedron.

This paper gives a solution to the problem of changing and
extending these intervals. We find the Feichtner-Kozlov concept of
building sets of arbitrary finite-meet semilattices, applied to
the case of simplicial complexes, appropriate for this purpose.

In particular, we are interested in iterating the procedure of
generating complexes of nested sets. The main obstacle is that
every iteration of this procedure increases the complexity of the
structure of nested curly braces in elements of the produced
simplicial complex. This makes it almost impossible to write down
the result of the procedure after just a couple of iterations. We
show here how to avoid this problem by using a simple algebraic
presentation of the members of nested sets as elements of a freely
generated commutative semigroup.

Such an iterated procedure is applied to the boundary complex of a
simplicial polytope. This provides an opportunity to replace the
simplex at the beginning of the interval simplex-permutohedron by
an arbitrary simple polytope, which automatically replaces the
permutohedron at the end of the interval by some other simple
polytope. By iterating the procedure, one obtains interesting
families of simplicial complexes.

Truncation of a polytope in its proper face is the operation that
guarantees that all these simplicial complexes may be realized as
the face lattices (with $\emptyset$ removed) of some simple
polytopes. In dimension 3, locally at each vertex, one may find
one of the following six types of truncation (cf.\ the examples of
Section 10).

\begin{center}
\begin{picture}(300,60)
\put(30,10){\line(0,1){30}} \put(30,40){\line(-2,1){30}}
\put(30,40){\line(2,1){30}}

\put(30,0){\makebox(0,0)[b]{\footnotesize no truncation}}

\put(150,10){\line(0,1){10}} \put(130,50){\line(-2,1){10}}
\put(170,50){\line(2,1){10}}

\put(150,20){\line(-2,3){20}} \put(150,20){\line(2,3){20}}
\put(130,50){\line(1,0){40}}

\put(150,0){\makebox(0,0)[b]{\footnotesize vertex}}

\put(265,10){\line(0,1){32.5}} \put(275,10){\line(0,1){32.5}}
\put(265,42.5){\line(-2,1){25}} \put(275,42.5){\line(2,1){25}}
\put(265,42.5){\line(1,0){10}}

\put(270,0){\makebox(0,0)[b]{\footnotesize edge}}

\end{picture}
\end{center}

\begin{center}
\begin{picture}(300,60)
\put(25,10){\line(0,1){17.5}} \put(35,10){\line(0,1){17.5}}
\put(10,50){\line(-2,1){10}} \put(50,50){\line(2,1){10}}

\put(25,27.5){\line(-2,3){15}} \put(35,27.5){\line(2,3){15}}
\put(10,50){\line(1,0){40}} \put(25,27.5){\line(1,0){10}}

\put(30,-1){\makebox(0,0)[b]{\footnotesize vertex + edge}}

\put(145,10){\line(0,1){17.5}} \put(155,10){\line(0,1){17.5}}
\put(130,50){\line(-2,1){10}} \put(155,50){\line(2,1){20}}
\put(165,42.5){\line(2,1){20}}

\put(145,27.5){\line(-2,3){15}} \put(155,27.5){\line(2,3){10}}
\put(130,50){\line(1,0){25}} \put(145,27.5){\line(1,0){10}}
\put(155,50){\line(4,-3){10}}

\put(150,-1){\makebox(0,0)[b]{\footnotesize vertex + 2 edges}}

\put(265,10){\line(0,1){17.5}} \put(275,10){\line(0,1){17.5}}
\put(265,50){\line(-2,1){20}} \put(255,42.5){\line(-2,1){20}}
\put(275,50){\line(2,1){20}} \put(285,42.5){\line(2,1){20}}

\put(265,27.5){\line(-2,3){10}} \put(275,27.5){\line(2,3){10}}
\put(265,50){\line(1,0){10}} \put(265,27.5){\line(1,0){10}}
\put(275,50){\line(4,-3){10}} \put(265,50){\line(-4,-3){10}}

\put(270,-1){\makebox(0,0)[b]{\footnotesize vertex + 3 edges}}

\end{picture}
\end{center}

\section{Nested complexes and abstract polytopes of hypergraphs}
In this section we compare the three approaches given in
\cite{FS05}, \cite{P09} and \cite{DP10c}, which all lead to the
same family of polytopes. We start with some preliminary notions.

Let $\alpha_1,\ldots,\alpha_n$ for $n\geq 1$ be finite and
incomparable (with respect to $\subseteq$) sets. Then
$C=P(\alpha_1)\cup\ldots\cup P(\alpha_n)$, where $P(\alpha)$ is
the power set of $\alpha$, is a \emph{finite abstract simplicial
complex} based on the set $\{\alpha_1,\ldots,\alpha_n\}$ of its
\emph{bases}\footnote{The term ``facet'' is reserved here for a
face of a polytope whose codimension is 1.}. Since we deal here
only with finite abstract simplicial complexes, we call them just
\emph{simplicial complexes}. For $C$ a simplicial complex, we have
that $\langle C,\cap\rangle$ is a meet-semilattice that induces
the poset $\langle C,\subseteq\rangle$. Depending on a context, by
a simplicial complex we usually mean such a meet-semilattice or a
poset.

For $B$ a subset of the domain of a function $f$, we denote by
$f[B]$ the set $\{f(b)\mid b\in B\}$. It is easy to see that two
simplicial complexes $C$ and $D$ are isomorphic as
meet-semilattices (or as posets) if and only if there is a
bijection $\varphi\!:\bigcup C\rightarrow\bigcup D$ such that
$\alpha\in C$ iff $\varphi[\alpha]\in D$ (or equivalently, such
that $\alpha$ is a basis of $C$ iff $\varphi[\alpha]$ is a basis
of $D$). If $\psi$ is a function such that the above $\varphi$ is
its restriction to $\bigcup C$, then we say that $\psi$
\emph{underlies} an isomorphism between $C$ and $D$. So, we have
the following remarks.

\prop{Remark 2.0.1}{If $C$ is a simplicial complex and $\psi$ is a
function such that its restriction to $\bigcup C$ is one-one, then
$\psi$ underlies an isomorphism between $C$ and the simplicial
complex $\{\psi[\alpha]\mid \alpha\in C\}$.}

\oprop{Remark 2.0.2}{If $\psi$ underlies isomorphisms between $C$
and $D_1$ and between $C$ and $D_2$, then $D_1=D_2$.}

Let $\alpha$ be an arbitrary finite set. According to Definition
7.1 of \cite{P09}, a set $B$ of nonempty subsets of $\alpha$ is a
\emph{building set}\footnote{Originally ``a building set on
$\alpha$'' in \cite{P09}, but I prefer to keep to the terminology
of \cite{FK04} where the building sets of semilattices are
introduced.} of $P(\alpha)$ when the following conditions hold:
\begin{tabbing}
\hspace{1.5em}\=(B1)\hspace{1em}\=If $\beta,\gamma\in B$ and
$\beta\cap \gamma\neq \emptyset$, then $\beta\cup \gamma\in B$;
\\[1ex]
\>(B2)\>$B$ contains all singletons $\{a\}$, for $a\in \alpha$.
\end{tabbing}
In the terminology of \cite{DP10c} (Sections 3-4), this notion
corresponds to the notion of atomic saturated hypergraph and it is
easy to see that this is a restriction to the case of power sets
of the more general notion of building sets of arbitrary
finite-meet semilattices given in \cite{FK04} (Definition 2.2).

Let $N$ be a family of sets. As in \cite{DP10c}, we say that
$\{\beta_1,\ldots, \beta_t\}\subseteq N$ is an
$N$-\emph{antichain} when $t\geq 2$ and $\beta_1,\ldots, \beta_t$
are incomparable with respect to $\subseteq$. Also, for a family
of sets $B$, we say that an $N$-antichain $\{\beta_1,\ldots,
\beta_t\}$ \emph{misses} $B$ when the union $\beta_1\cup\ldots\cup
\beta_t$ does not belong to $B$.

\subsection{Nested set complexes of Feichtner and Sturmfels}
Let $B$ be a building set of $P(\alpha)$, containing $\alpha$. A
subset $N$ of $B$ is \emph{nested} when every $N$-antichain
\emph{misses} B. It is easy to see that for $M\subseteq N\subseteq
B$, if $N$ is nested, then $M$ is nested too. So, the nested
subsets of $B$ form a simplicial complex whose bases are the
maximal nested subsets of $B$.

In \cite{FS05}, for $B$ as above, this simplicial complex would be
denoted by $\Nt(P(\alpha),B)$, and the link of $\alpha$ in
$\Nt(P(\alpha),B)$ would be denoted by $\N(P(\alpha),B)$ and
called the \emph{nested set complex} of $P(\alpha)$ with respect
to $B$. (Nested set complexes of \cite{FS05} are defined not only
for power sets but for arbitrary finite lattices.)

Unfortunately, the name nested set complex and the symbol $\N$
were previously used in \cite{FM03} (see also \cite{C07}) for the
result and the name of the operation that generalizes $\Nt$ from
above.

\subsection{Nested complexes of Postnikov}
Let $B$ be a building set of $P(\alpha)$ not necessarily
containing $\alpha$. According to Definition 7.3 of \cite{P09}, a
subset $N$ of $B$ is a \emph{nested set} when it satisfies the
following:
\begin{tabbing}
\hspace{1.5em}\=(N1)\hspace{1em}\=If $\beta,\gamma\in N$ then
($\beta\subseteq \gamma$ or $\gamma\subseteq\beta$ or
$\beta\cap\gamma=\emptyset$);
\\[1ex]
\>(N2)\>If $\beta_1,\ldots, \beta_t$ for $t\geq 2$ are mutually
disjoint elements of $N$, then
\\
\>\>$\beta_1\cup\ldots\cup \beta_t$ does not belong to $B$;
\\[1ex]
\>(N3)\>$N$ contains all maximal elements of $B$.
\end{tabbing}

The \emph{nested complex} $\N(B)$ is a poset of all nested sets
ordered by inclusion. Let $\N^\ast(B)$ be obtained from $\N(B)$ by
removing every maximal element of $B$ from nested sets. Then
$\N^\ast(B)$ is a simplicial complex whose bases are maximal
nested sets with maximal elements of $B$ removed.

\subsection{Abstract polytopes of hypergraphs}

This is an alternative, inductive approach to the same matters,
which is given in \cite{DP10c}. For $\alpha$ a finite set, let
$H\subseteq P(\alpha)$ be such that $\emptyset\not\in H$ and
$\alpha=\bigcup H$. Then $H$ is a \emph{hypergraph} on $\alpha$
(see \cite{B89}, Section 1.1). A hypergraph $H$ is \emph{atomic}
when for every $x\in \bigcup H$ we have that $\{x\}\in H$ (see
\cite{DP10c}, Section~3).

A \emph{hypergraph partition} of a hypergraph $H$ is a partition
$\{H_1,\ldots,H_n\}$, with $n\geq 1$, of $H$ such that $\{\bigcup
H_1,\ldots,\bigcup H_n\}$ is a partition of $\bigcup H$. A
hypergraph $H$ is \emph{connected} when it has only one hypergraph
partition. (If $H$ is nonempty, then this hypergraph partition is
$\{H\}$.) For example, $H=\{\{x\},\{y\},\{z\},\{x,y,z\}\}$ is
connected.

A hypergraph partition $\{H_1,\ldots, H_n\}$ of $H$ is
\emph{finest} when every $H_i$ is a connected hypergraph on
$\bigcup H_i$. We say that $\bigcup H_i$ is a \emph{connected
component} of $\bigcup H$. For $H\subseteq P(\alpha)$ and
$\beta\subseteq\alpha$ let
\[
H_\beta=_{df}\{\gamma\in H\mid \gamma\subseteq\beta\}=H\cap
P(\beta).
\]

Let $H$ be an atomic hypergraph. By induction on the cardinality
$|\bigcup H|$ of $\bigcup H$ the \emph{constructions} of $H$ are
defined as follows
\begin{itemize}
\item[(0)] if $|\bigcup H|=0$, then $H$ is the empty hypergraph
$\emptyset$, and $\emptyset$ is the only construction of
$\emptyset$;

\item[(1)] if $|\bigcup H|\geq 1$, and $H$ is connected, and $K$
is a construction of $H_{\ccup H-\{x\}}$ for $x\in\bigcup H$, then
$K\cup\{\bigcup H\}$ is a construction of~$H$;

\item[(2)] if $|\bigcup H|\geq 2$, and $H$ is not connected, and
$\{H_1,\ldots,H_n\}$, where $n\geq 2$, is the finest hypergraph
partition of $H$, and for every $i\in\{1,\ldots,n\}$ we have that
$K_i$ is a construction of $H_i$, then $K_1\cup\ldots\cup K_n$ is
a construction of~$H$.
\end{itemize}

A subset of a construction is called a \emph{construct} when it
contains every connected component of $\bigcup H$. The abstract
polytope of $H$, denoted by $\A(H)$, is a poset of all the
constructs of $H$ ordered by $\supseteq$, plus the bottom element.

A hypergraph $H$ is \emph{saturated} when it satisfies the
condition (B1) of the definition of building set with $B$ replaced
by $H$. So, the notions of atomic saturated hypergraph on $\alpha$
and of building set of $P(\alpha)$ coincide. It follows from
Proposition 4.7 of \cite{DP10c} that for $\bar{H}$ being the
intersection of all the atomic saturated hypergraphs, i.e.\
building sets, containing an atomic hypergraph $H$, we have
$\A(H)=\A(\bar{H})$.

\subsection{The interval simplex-permutohedron}

Our task now is to compare the posets $\N(P(\alpha),B)$, $\N(B)$
and $\A(B)$ for $B$ a building set of $P(\alpha)$. We start with
some abbreviations. Let $L_-$, for a finite lattice $L$ be the
poset obtained from $L$ by removing the bottom element and,
analogously, let $L^-$ be obtained from $L$ by removing the top
element. Let $P^{op}$, for a poset $P$ be the poset obtained by
reversing the order of $P$.

Let $B$ be a building set of $P(\alpha)$ and let $N\subseteq B$.
Then it is easy to verify that
\begin{tabbing}
\hspace{1.5em}($\ast$)\quad (N1) and (N2) of Section 2.2 hold iff
every $N$-antichain misses $B$.
\end{tabbing}
From ($\ast$) we infer that for $B$ containing $\alpha$, we have
the isomorphism between $\N(B)$ and $\N(P(\alpha),B)$ obtained by
removing $\alpha$ from the elements of $\N(B)$.

Although Proposition 6.12 of \cite{DP10c} is formulated for
atomic, saturated and connected hypergraphs, it is easy to verify
that connectedness is not essential, i.e.\ that the following
proposition holds.

\prop{Proposition 2.4.1}{For every building set $H$ and every
$N\subseteq H$ we have: every $N$-antichain misses $H$ iff for
some construction $K$ of $H$ we have that $N\subseteq K$.}

\noindent From Proposition 2.4.1 and ($\ast$) we infer that
$(\N(B))^{op}$ is equal to $(\A(B))_-$.

Since $\N(B)$ and $\A(B)$ for $B$ not containing $\alpha$ are
easily defined in terms of $\N(H_1),\ldots,\N(H_n)$, for
$\{H_1,\ldots,H_n\}$ being the finest hypergraph partition of $B$
(see the definition of $\otimes$ given in \cite{DP10c}, Section~5)
we may consider only the building sets of $P(\alpha)$ that contain
$\alpha$ (if $\alpha$ is nonempty). Such a building set is called
in \cite{DP10c} atomic saturated connected hypergraph, or shortly
ASC-hypergraph. For $H$ being an ASC-hypergraph, all the three
posets $(\N(P(\alpha),H))^{op}$, $(\N(H))^{op}$ and $(\A(H))_-$
are isomorphic.

Let $\alpha$ be finite, nonempty set. The ASC-hypergraphs on
$\alpha$, i.e.\ the building sets of $P(\alpha)$ containing
$\alpha$, make a meet-semilattice $\H$ with $\cap$ being the meet
operation. We have that $H_\bot=\{\{a\}\mid
a\in\alpha\}\cup\{\alpha\}$ and $H_\top=P(\alpha)-\{\emptyset\}$
are respectively the bottom and the top element of $\H$. (Since
every pair of elements of $\H$ has the upper bound, it is easy to
define the lattice structure on $\H$.)

Let $n=|\alpha|-1$. We have that $\A(H_\bot)$ is isomorphic to the
face lattice of an $n$-dimensional simplex. On the other hand,
$\A(H_\top)$ is isomorphic to the face lattice of an
$n$-dimensional permutohedron. For every $H\in\H$ we have that
$\A(H)$ is isomorphic to the face lattice of some $n$-dimensional
polytope $\G(H)$ (see \cite{DP10c}, Section~9 and Appendix~B). For
a fixed dimension $n\geq 0$, we define the interval
\emph{simplex-permutohedron} to be the set $\{\G(H)\mid H\in\H\}$.

Since the function $\G$ defined in \cite{DP10c}, Section~9, is
one-one (although some members of the family are combinatorially
equivalent), the interval simplex-permutohedron may be enriched
with a poset or a lattice structure induced by the structure of
$\H$. However, if we consider, as usual, polytopes only up to
combinatorial equivalence, then the interval simplex-permutohedron
is just a family of polytopes. We refer to \cite{DP10c},
appendix~B, where one may find the complete intervals
simplex-permutohedron in dimension 3 and lower, and a chart of a
fragment of $\H$ tied to dimension 3.

\section{Complexes of nested sets for simplicial complexes}
In this section we define the notion of a building set of a
simplicial complex and the notion of a corresponding complex of
nested sets, which are just restrictions of the notions defined in
\cite{FK04} for arbitrary finite-meet semilattices.

Let $C$ be a simplicial complex based on the set
$\{\alpha_1,\ldots,\alpha_n\}$ (cf.\ Section~2). A set $B\subseteq
C$ is a \emph{building set} of $C$ when for every
$i\in\{1,\ldots,n\}$, we have that $B_{\alpha_i}=B\cap
P(\alpha_i)$ is a building set of $P(\alpha_i)$ in the sense of
the definition given at the beginning of Section~2. It is easy to
verify that this is the notion to which the notion of a building
set of a finite-meet semilattice introduced in \cite{FK04}
(Definition 2.2) is reduced to, when the semilattice is a
simplicial complex.

Let $\min X$, for $X$ a family of sets, be the set of minimal
(with respect to $\subseteq$) members of $X$, and let $\max X$ be
defined analogously. Then we can prove the following proposition
in a straightforward manner.

\prop{Proposition 3.1}{For a finite, nonempty set $\alpha$ and a
simplicial complex $C$, we have that:
\\[1ex]
{\rm (1)}\quad $B$ is a building set of $C$ iff $B_\gamma$ is a
building set of $P(\gamma)$ for every $\gamma\in C$;
\\[1ex]
{\rm (2)}\quad if $B$ is a building set of $P(\alpha)$, then
$B-\{\alpha\}$ is a building set of $P(\alpha)-\{\alpha\}$;
\\[1ex]
{\rm (3)}\quad $B_\bot=\{\{a\}\mid a\in \bigcup C\}$ is a building
set of $C$;
\\[1ex]
{\rm (4)}\quad if $\alpha\in C$, then $B_\bot\cup\{\alpha\}$ is a
building set of $C$;
\\[1ex]
{\rm (5)}\quad if $B$ is a building set of $C$ and
$\beta\in\min(B-B_\bot)$, then $B-\{\beta\}$ is a building set of
$C$;
\\[1ex]
{\rm (6)}\quad if $B$ is a building set of $C$ and
$\beta\in\max(C-B)$, then $B\cup\{\beta\}$ is a building set of
$C$. }

\noindent (Note that if $\alpha$ is finite and nonempty, then
$P(\alpha)-\{\alpha\}$ is a simplicial complex based on the set
$\{\alpha-\{a\}\mid a\in\alpha\}$.)

Let $B$ be a building set of a simplicial complex $C$. A subset
$N$ of $B$ is \emph{nested} when for every $N$-antichain
$\{\beta_1,\ldots, \beta_t\}$, the union $\beta_1\cup\ldots\cup
\beta_t$ belongs to $C-B$. This is again in accordance with
Definition 2.7 of \cite{FK04}.

It is easy to see that for $M\subseteq N\subseteq B$, if $N$ is
nested, then $M$ is nested too. So, the nested subsets of $B$ form
again a simplicial complex whose bases are the maximal nested
subsets of $B$. We denote this simplicial complex by $\Nt(C,B)$.
(According to \cite{FK04} it would be just $\N(B)$, but since we
want to make a comparison with the notions given in Section 2.1,
we find $\Nt(C,B)$ more appropriate since it is reduced to
$\Nt(P(\alpha),B)$ when $C$ is $P(\alpha)$ and $B$ contains
$\alpha$.) We have that $\bigcup\Nt(C,B)=B$, since for every
$\beta\in B$ we have that $\{\beta\}$ is nested. Then we can prove
the following proposition in a straightforward manner.

\prop{Proposition 3.2}{For $B$ a building set of a simplicial
complex $C$ the following statements are equivalent:
\\[1ex]
{\rm (1)}\quad $N\in\Nt(C,B)$ ;
\\[1ex]
{\rm (2)}\quad $N\subseteq B$, $\bigcup N\in C$ and every
$N$-antichain misses $B$;
\\[1ex]
{\rm (3)}\quad there is $\gamma\in C$ such that
$N\in\Nt(P(\gamma),B_\gamma)$;
\\[1ex]
{\rm (4)}\quad there is a basis $\alpha$ of $C$ such that
$N\in\Nt(P(\alpha),B_\alpha)$. }

\noindent (Note that for (3) and (4) above, we do not require that
$B_\gamma$ contains $\gamma$ and that $B_\alpha$ contains $\alpha$
as in the definition of $\Nt$ given in Section 2.1.)

The following proposition is a corollary of Propositions 2.4.1 and
3.2.

\prop{Proposition 3.3}{We have that $N\in\Nt(C,B)$ iff there is
$\gamma\in C$ such that $N$ is a subset of a construction of the
hypergraph $B_\gamma$. Moreover, $N\in\Nt(C,B)$ is a basis of
$\Nt(C,B)$ iff there is a basis $\alpha$ of $C$ such that $N$ is a
construction of the hypergraph $B_\alpha$.}

The following proposition defines $\N$ of Section 2.1 in terms of
$\Nt$ from above.

\propo{Proposition 3.4}{If $B$ is a building set of $P(\alpha)$,
containing $\alpha$, then
\[
\N(P(\alpha),B)=\Nt(P(\alpha)-\{\alpha\},B-\{\alpha\}). \] }

\dkz Since $\alpha\in B$ we have that $\alpha$ is nonempty. By (2)
of Proposition 3.1, $B-\{\alpha\}$ is a building set of the
simplicial complex $P(\alpha)-\{\alpha\}$.

For the proof of $\subseteq$-direction, let $N\in\N(P(\alpha),B)$.
By the definition of $\N$ (cf.\ Section 2.1), we have
$N\in\Nt(P(\alpha),B)$ and $N$ is in the link of $\alpha$. Hence,
$N\subseteq B-\{\alpha\}$. We have that $\bigcup N$ is either
empty, or it is equal to $\beta$ for some $\beta\in N$, or it is
equal to $\beta_1\cup\ldots\cup\beta_t$ for some $N$-antichain
$\{\beta_1,\ldots,\beta_t\}$. In all the three cases it follows
that $\bigcup N\in P(\alpha)-\{\alpha\}$, since $\alpha\not\in N$,
$\alpha\in B$ and $N\in\Nt(P(\alpha),B)$. Let
$\{\beta_1,\ldots,\beta_t\}$ be an $N$-antichain. From
$N\in\Nt(P(\alpha),B)$ we conclude that
$\beta_1\cup\ldots\cup\beta_t\not\in B$, and hence
$\beta_1\cup\ldots\cup\beta_t\not\in B-\{\alpha\}$. So, by (2) of
Proposition 3.2, we have that
$N\in\Nt(P(\alpha)-\{\alpha\},B-\{\alpha\})$.

For the proof of $\supseteq$-direction, let
$N\in\Nt(P(\alpha)-\{\alpha\},B-\{\alpha\})$. So, $N\subseteq
B-\{\alpha\}$ and hence $N\subseteq B$. If
$\{\beta_1,\ldots,\beta_t\}$ is an $N$-antichain, then
\[
\beta_1\cup\ldots\cup\beta_t\in(P(\alpha)-\{\alpha\})-
(B-\{\alpha\})=P(\alpha)-B.
\]
Hence, $N\in\Nt(P(\alpha),B)$ and analogously, since every
$(N\cup\{\alpha\})$-antichain is an $N$-antichain, we have that
$N\cup\{\alpha\}\in\Nt(P(\alpha),B)$. From this and the fact that
$\alpha\not\in N$, we conclude that $N$ belongs to the link of
$\alpha$ in $\Nt(P(\alpha),B)$. Hence $N\in\N(P(\alpha),B)$. \qed

This proposition sheds a new light on the interval
simplex-permutohedron and it points out a way how to change the
boundaries of this interval. We will discuss this later in
Section~10.

\section{Commutative semigroup notation}
We show how to use the elements of commutative semigroups freely
generated by sets to interpret the families of families...of sets.
We deal with ``formal sums'' of elements of some set and with
functions mapping finite, nonempty subsets of this set to the
formal sums, or mapping families of finite, nonempty subsets to
families of formal sums etc. For example, if we take the set
$\{x,y,z,u\}$, then we are interested in mappings of the following
form:
\begin{tabbing}
\hspace{2em}\=$\{x,y,z\}\mapsto x+y+z$,
\\[.5ex]
\>$\{\{x\},\{x,y\},\{z\}\}\mapsto\{x,x+y,z\}$, and
\\[.5ex]
\>$\{\{\{x\}\},\{\{x\},\{y,u\}\}\}\mapsto\{\{x\},\{x,y+u\}\}$.
\end{tabbing}

For the sake of precision, we use some basic notions from category
theory, which all may be found in \cite{ML98}. However, no result
from category theory is needed in the sequel. We try to illustrate
all the notions we use by examples so that a reader not familiar
with category theory may deal with these notions intuitively.

Let $\textbf{Cs}$ be the category of commutative semigroups and
let $G\!:\textbf{Cs}\rightarrow \textbf{Set}$ be the forgetful
functor, which maps a commutative semigroup to its underlying set.
Let $F\!:\textbf{Set}\rightarrow\textbf{Cs}$ be its left adjoint
assigning to every set the free commutative semigroup generated by
this set. (If $X$ is empty, then $FX$ is the empty semigroup.)
Consider the monad $\langle T,\eta,\mu\rangle$ defined by this
adjunction ($T=G\cirk F$, $\eta$ is the unit of this adjunction
and $\mu=G\varepsilon F$, where $\varepsilon$ is the counit of
this adjunction).

For every set $X$ we have that $TX$ is the underlying set of the
free commutative semigroup generated by $X$. We can take the
elements of $TX$ to be the formal sums of elements of $X$, so that
if $X=\{x,y,z,u\}$, then
\[
TX=\{x,y,z,u,2x,x+y,\ldots,3x+y+5z+u,\ldots\}.
\]
We use $+^2$ to denote the formal addition of $T^2X=TTX$.
Analogously, $+^3$ denotes the formal addition of $T^3X$ etc. So,
if $X$ is as above, then
\[
T^2X=\{x,\ldots,3x+y+5z+u,\ldots, x+^2
2x,\ldots,z+^2(3x+y+5z+u),\ldots\}.
\]
For the natural transformations $\eta$ and $\mu$ of the monad
$\langle T,\eta,\mu\rangle$, we have that $\eta_X:X\rightarrow TX$
is such that $\eta_X(x)=x$, for every $x\in X$, and
$\mu_X:T^2X\rightarrow TX$ evaluates $+^2$ as $+$, so that, for
example, $\mu_X(z+^2(3x+y+5z+u))=3x+y+6z+u$.

Besides the endofunctor $T\!:\textbf{Set}\rightarrow\textbf{Set}$
we are interested in the power set functor
$P\!:\textbf{Set}\rightarrow\textbf{Set}$ (see \cite{ML98},
Section I.3) and its modification
$P_F\!:\textbf{Set}\rightarrow\textbf{Set}$, which assigns to a
set $X$ the set of finite, nonempty subsets of $X$. On arrows,
$P_F$ is defined as the power set functor, so that for
$f\!:X\rightarrow Y$ and $A\in P_FX$, we have that $P_Ff(A)=f[A]$.

For every set $X$, let $\sigma_X\!:P_FX\rightarrow TX$ be the
function such that for $A=\{a_1,\ldots,a_n\}\subseteq X$, where
$n\geq 1$, we have $\sigma_X(A)=a_1+\ldots + a_n$. It follows that
$\sigma_X$ is one-one. (Note that $\sigma$ is not a natural
transformation from $P_F$ to $T$.) For example, if $X$ is as
above, then $\sigma_X(\{x\})=x$, and for
$N=\{\{x\},\{x,y\},\{z\}\}$, we have that
$P\sigma_X(N)=\{x,x+y,z\}$.

\section{Iterating $\Nt$}
If $B_0$ is a building set of a simplicial complex $C$, then
$\Nt(C,B_0)$ is a simplicial complex, and if $B_1$ is a building
set of $\Nt(C,B_0)$, then $\Nt(\Nt(C,B_0),B_1)$ is again a
simplicial complex and we may iterate this procedure. So, for
$n\geq 0$, $B_0$ a building set of $C$, and $B_{i+1}$, for $0\leq
i\leq n$ a building set of $\Nt(C,B_0,\ldots,B_i)$, we define
\[
\Nt(C,B_0,\ldots,B_{n+1})=_{df}\Nt(\Nt(C,B_0,\ldots,B_n),B_{n+1}).
\]

\noindent {\sc Example} 5.1\quad Let $C$ be a simplicial complex
whose bases are $\{x,y,z\}$, $\{x,y,u\}$, $\{x,z,u\}$ and
$\{y,z,u\}$, i.e.\ $C=P(X)-\{X\}$ for $X=\{x,y,z,u\}$, and let
$B_0=\{\{x\},\{y\},\{z\},\{u\},\{x,y\}\}$. By Proposition 3.3, the
bases of $\Nt(C,B_0)$ are the constructions of $(B_0)_\alpha$ for
all the bases $\alpha$ of $C$. So, $\Nt(C,B_0)$ has the following
six bases
\[
\begin{array}{lll}
\{\{x\},\{x,y\},\{z\}\} \;\;\;\mbox{\rm and} &
\{\{y\},\{x,y\},\{z\}\} & \mbox{\rm
derived from }\{x,y,z\},\\[.5ex]
\{\{x\},\{x,y\},\{u\}\} \;\;\;\mbox{\rm and} &
\{\{y\},\{x,y\},\{u\}\} & \mbox{\rm
derived from }\{x,y,u\},\\[.5ex]
& \{\{x\},\{z\},\{u\}\} & \mbox{\rm derived from }\{x,z,u\},
\\[.5ex]
& \{\{y\},\{z\},\{u\}\} & \mbox{\rm derived from }\{y,z,u\}.
\end{array}
\]
For $B_1=\{\{\{x\}\},\{\{y\}\},\{\{z\}\},\{\{u\}\},
\{\{x,y\}\},\{\{x\},\{x,y\}\}\}$, again by Proposition 3.3, we
have that $\Nt(C,B_0,B_1)$ has the following eight bases.
\[\begin{array}{ll}
\{\{\{x\}\},\{\{x\},\{x,y\}\},\{\{z\}\}\} &
\{\{\{x,y\}\},\{\{x\},\{x,y\}\},\{\{z\}\}\}
\\
\{\{\{x\}\},\{\{x\},\{x,y\}\},\{\{u\}\}\} &
\{\{\{x,y\}\},\{\{x\},\{x,y\}\},\{\{u\}\}\}
\\
\{\{\{y\}\},\{\{x,y\}\},\{\{z\}\}\} &
\{\{\{y\}\},\{\{x,y\}\},\{\{u\}\}\}
\\
\{\{\{x\}\},\{\{z\}\},\{\{u\}\}\} &
\{\{\{y\}\},\{\{z\}\},\{\{u\}\}\}
\end{array}
\]

This example illustrates how the iterated application of $\Nt$
increases complexity of the structure of nested curly braces in
the elements of the resulting simplicial complex. Our intention is
to make these matters simpler.

For $C$ a simplicial complex such that $\bigcup C\subseteq X$, $B$
a building set of $C$, and $\sigma_X\!:P_FX\rightarrow TX$ defined
in the preceding section, let
\[
\St_X(C,B)=_{df}\{P\sigma_X(N)\mid N\in\Nt(C,B)\}.
\]
It is easy to see that $\St_X(C,B)$ is indeed a simplicial
complex. For $C$, $X$ and $B_0$ as in Example 5.1, we have that
$\St_X(C,B_0)$ is a simplicial complex with the bases
$\{x,x+y,z\}$, $\{y,x+y,z\}$, $\{x,x+y,u\}$, $\{y,x+y,u\}$,
$\{x,z,u\}$ and $\{y,z,u\}$. We have the following proposition.

\prop{Proposition 5.2}{The function $\sigma_X$ underlies an
isomorphism between $\Nt(C,B)$ and $\St_X(C,B)$.}

\dkz Since $\sigma_X$ is one-one, we just apply Remark 2.0.1. \qed

Let $\sigma_X^1=_{df}\sigma_X$, and for $n\geq 1$, let
$\sigma_X^{n+1}=_{df}\sigma_{T^nX}\cirk P_F\sigma_X^n$. It is easy
to see that $\sigma_X^n\!:P_F^nX\rightarrow T^nX$ satisfies the
following inductive clauses:
\begin{tabbing}
\hspace{1.5em}\=$\sigma_X^1(\{a_1,\ldots,a_k\})=a_1+\ldots+a_k$,
\\[1ex]
\>$\sigma_X^{n+1}(\{\alpha_1,\ldots,\alpha_m\})=\sigma_X^n(\alpha_1)+^{n+1}\ldots
+^{n+1}\sigma_X^n(\alpha_m)$.
\end{tabbing}
For example, we have that $\sigma_X^2(\{\{x\}\})=x$,
$\sigma_X^2(\{\{x\},\{x,y\},\})=x+^2(x+y)$, and
$\sigma_X^3(\{\{\{x\}\},\{\{x\},\{x,y\}\}\})=x+^3(x+^2(x+y))$.

For $n\geq 1$, $C$ a simplicial complex such that $\bigcup
C\subseteq X$, $B_0$ a building set of $C$, and $B_i$, for $1\leq
i\leq n$ a building set of $\Nt(C,B_0,\ldots,B_{i-1})$, let
\[
\St_X(C,B_0,\ldots,B_n)=_{df}\{P\sigma_X^{n+1}(N)\mid
N\in\Nt(C,B_0,\ldots,B_n)\}.
\]

We easily compute that $\St_X(C,B_0,B_1)$, for $C$, $X$, $B_0$ and
$B_1$ as in Example 5.1, has the following eight bases
$\{x,x+^2(x+y),z\}$, $\{x+y,x+^2(x+y),z\}$, $\{x,x+^2(x+y),u\}$,
$\{x+y,x+^2(x+y),u\}$, $\{y,x+y,z\}$, $\{y,x+y,u\}$, $\{x,z,u\}$,
$\{y,z,u\}$. The following proposition is a generalization of
Proposition 5.2.

\prop{Proposition 5.3}{The function $\sigma_X^{n+1}$ underlies an
isomorphism between $\Nt(C,B_0,\ldots,B_n)$ and
$\St_X(C,B_0,\ldots,B_n)$.}

\dkz Since $\sigma_X$ is one-one and $P_F$ preserves this
property, we conclude that $\sigma_X^{n+1}$ is one-one. Then we
apply Remark 2.0.1. \qed

\propo{Remark 5.4}{If $\psi$ is a function that underlies an
isomorphism between simplicial complexes $C$ and $D$, and $B$ is a
building set of $C$, then $P^2\psi(B)$ is a building set of $D$
and the function $P_F\psi$ underlies an isomorphism between
$\Nt(C,B)$ and $\Nt(D,P^2\psi(B))$.}

\prop{Proposition 5.5}{For $n\geq 1$, the function
$\sigma_X^{n+1}$ underlies an isomorphism between
$\Nt(C,B_0,\ldots,B_n)$ and
$\St_{T^nX}(\St_X(C,B_0,\ldots,B_{n-1}),P^2\sigma_X^n(B_n))$.}

\dkz By Proposition 5.3 and Remark 5.4 we have that
$P_F\sigma_X^n$ underlies an isomorphism between
$\Nt(\Nt(C,B_0,\ldots,B_{n-1}),B_n)$, which is equal to
$\Nt(C,B_0,\ldots,B_n)$ and
$\Nt(\St_X(C,B_0,\ldots,B_{n-1}),P^2\sigma_X^n(B_n))$. By
Proposition 5.2, since
$\bigcup\St_X(C,B_0,\ldots,B_{n-1})\subseteq T^nX$, we have that
$\sigma_{T^nX}$ underlies an isomorphism between
$\Nt(\St_X(C,B_0,\ldots,B_{n-1}),P^2\sigma_X^n(B_n))$ and
\linebreak
$\St_{T^nX}(\St_X(C,B_0,\ldots,B_{n-1}),P^2\sigma_X^n(B_n))$. So,
$\sigma_{T^nX}\cirk P_F\sigma_X^n$, which is equal to
$\sigma_X^{n+1}$ underlies the desired isomorphism. \qed

For $n\geq 1$, by Propositions 5.3, 5.5 and Remark 2.0.2 we have
that
\[\St_X(C,B_0,\ldots,B_n)=\St_{T^nX}(\St_X(C,B_0,\ldots,B_{n-1}),P^2\sigma_X^n(B_n)).\]

Proposition 5.3 enables us to simplify a bit the notation---we
just codify complicated curly braces notation using $+$, $+^2$,
$+^3$, etc., leaving the original nested structure untouched.
However, we will show that the necessity of such a nested
structure for describing these simplicial complexes is an
illusion.

Let again $C$, $X$, $B_0$ and $B_1$ be as in Example 5.1. It is
easy to check that $\St_X(C,B_0,B_1)$ (and hence $\Nt(C,B_0,B_1)$)
is isomorphic to the simplicial complex $K$ whose bases are
$\{x,2x+y,z\}$, $\{x+y,2x+y,z\}$, $\{x,2x+y,u\}$,
$\{x+y,2x+y,u\}$, $\{y,x+y,z\}$, $\{y,x+y,u\}$, $\{x,z,u\}$ and
$\{y,z,u\}$. (To obtain $K$ from $\St_X(C,B_0,B_1)$, we have just
``evaluated'' $+^2$ as $+$ in the elements of its bases.) We show
that this will always be the case.

\section{Flattening the nested structure}
For a component $\mu_X:T^2X\rightarrow TX$ of the multiplication
of the monad $\langle T,\eta,\mu\rangle$ (cf.\ Section~4), we
define inductively the function $\mu_X^n\!:T^{n+1}X\rightarrow TX$
in the following way. For $n=0$, let $\mu_X^0\!:TX\rightarrow TX$
be the identity function on $TX$, and for $n\geq 1$, let
$\mu_X^n=_{df} \mu_X\cirk T\mu_X^{n-1}$. (For $n\geq 1$, by the
associativity law of the monad, we have that $\mu_X^n=\mu_X\cirk
T\mu_X\cirk\ldots\cirk
T^{n-1}\mu_X=\mu_X\cirk\mu_{TX}\cirk\ldots\cirk\mu_{T^{n-1}X}$.)
For example, $\mu_X^1(x+^2(x+y))=2x+y$ and
$\mu_X^2(x+^3(x+^2(x+y)))=3x+y$.

For $\St_X(C,B_0,\ldots,B_n)$ defined in the preceding section, we
define
\[
\Ss_X(C,B_0,\ldots,B_n)=_{df}\{P\mu_X^n(\alpha)\mid
\alpha\in\St_X(C,B_0,\ldots,B_n)\}.
\]
In our example above, we have that $K=\Ss_X(C,B_0,B_1)$ for $C$,
$X$, $B_0$ and $B_1$ as in Example 5.1. Since the function
$\mu_X^n$ is not one-one when $X$ is nonempty and $n\geq 1$, it is
not so obvious that $\mu_X^n$ underlies an isomorphism between
$\St_X(C,B_0,\ldots,B_n)$ and $\Ss_X(C,B_0,\ldots,B_n)$. The rest
of this section is devoted to a proof of this result.

For $C$ a simplicial complex such that $\bigcup C\subseteq TX$ and
$B$ a building set of $C$, we formulate some conditions, which
guarantee that $\mu_X$ restricted to $\bigcup\St_{TX}(C,B)$ is
one-one. That this is not always the case is shown by the
following example.

Let $X=\{x,y\}$ and $C=P(\{x,x+y,2x+y\})$. For a building set
$B=\{\{x\},\{x+y\},\{2x+y\},\{x,x+y\}\}$ of $C$, we have that
$\bigcup\St_{TX}(C,B)=\{x,x+y,2x+y,x+^2(x+y)\}$ and
$\mu_X(2x+y)=\mu_X(x+^2(x+y))=2x+y$.

Let $f$ be a function, and let $C$ be a simplicial complex such
that for every $a\in\bigcup C$ we have that $f(a)\in
\textbf{R}^n$. We say that $f$ \emph{faithfully realizes} $C$ when
for every $\alpha,\beta\in C$, every $\{k_a>0\mid a\in \alpha$\}
and every $\{l_b>0\mid b\in\beta\}$, we have that
\[
\sum_{a\in\alpha}k_a\cdot f(a)=\sum_{b\in\beta}l_b\cdot f(b),
\]
implies that $\alpha=\beta$ and $k_a=l_a$ for every $a\in\alpha$.

For $f$ and $C$ as above, let
\[
f(C)=_{df}\{\mbox{\rm
cone}(\{f(a)\mid a\in\alpha\})\mid\alpha\in C\}.
\]
(See \cite{Z95} for the notions of cone and simplicial fan.) It is
easy to conclude that if $f$ faithfully realizes $C$, then $f(C)$
is a simplicial fan whose face lattice is isomorphic to $C$. On
the other hand, if $f(C)$ is a simplicial fan whose face lattice
is isomorphic to $C$ via the mapping that extends $\{a\}\mapsto
\mbox{\rm cone}(\{f(a)\})$ for $a\in\bigcup C$, then $f$
faithfully realizes $C$. We can now prove the following.

\prop{Proposition 6.1}{Suppose $f\!:A\rightarrow \textbf{R}^n$ and
$g\!:B\rightarrow \textbf{R}^n$ faithfully realize $C$ and $D$
respectively, and let $\bigcup g(D)\subseteq \bigcup f(C)$. If
$L\!:\textbf{R}^n\rightarrow \textbf{R}^m$ is a linear
transformation such that $Lf$ faithfully realizes $C$, then $Lg$
faithfully realizes~$D$. }

\dkz Suppose $\gamma,\delta\in D$, $p_c,q_d>0$ for $c\in\gamma$,
$d\in\delta$, and
\[
\sum_{c\in\gamma}p_c\cdot Lg(c)=\sum_{d\in\delta}q_d\cdot Lg(d).
\]
For $x=\sum_{c\in\gamma}p_c\cdot g(c)$ and
$y=\sum_{d\in\delta}q_d\cdot g(d)$, we have $x,y\in\bigcup
g(D)\subseteq \bigcup f(C)$. So, there are $\alpha,\beta\in C$,
$k_a,l_b>0$ for $a\in\alpha$ and $b\in\beta$, such that
$x=\sum_{a\in\alpha}k_a\cdot f(a)$ and $y=\sum_{b\in\beta}l_b\cdot
f(b)$.

From $L(x)=L(y)$ we obtain
\[
\sum_{a\in\alpha}k_a\cdot Lf(a)=\sum_{b\in\beta}l_b\cdot Lf(b),
\]
and since $Lf$ faithfully realizes $C$, we have $\alpha=\beta$ and
$k_a=l_a$ for every $a\in\alpha$. Hence, $x=y$ and since $g$
faithfully realizes $D$ we obtain $\gamma=\delta$ and $p_c=q_c$
for every $c\in\gamma$. \qed

Let $X$ be a finite and nonempty set. We always assume that $X$ is
linearly ordered. In our examples, if $X=\{x,y,z,u\}$ then we
assume $x<y<z<u$. There is a natural identification of $TX$ with
$\textbf{N}^{|X|}-\{(0,\ldots,0)\}$ when $X$ is finite and
linearly ordered. If $X$ is as above, then we identify $2x$ with
$(2,0,0,0)$, $3x+y+5z+u$ with $(3,1,5,1)$ etc. Let
$\kappa_X\!:TX\rightarrow \textbf{R}^{|X|}$ be defined by this
identification.

\prop{Lemma 6.2}{Let $A,B\subseteq TX$ be such that $\sigma_X
(X)\in A\cap B$ and that there are $x,y\in X$ (not necessarily
distinct) such that $A-\{\sigma_X(X)\}\subseteq T(X-\{x\})$ and
$B-\{\sigma_X(X)\}\subseteq T(X-\{y\})$. If
\[
(\ast)\quad \sum_{a\in A}k_a\cdot \kappa_X(a)=\sum_{b\in
B}l_b\cdot\kappa_X(b),\;\mbox{\rm for }k_a,l_b>0,
\]
then $A-\{\sigma_X(X)\},B-\{\sigma_X(X)\}\subseteq
T(X-(\{x\}\cup\{y\}))$ and $k_{\sigma_X(X)}=l_{\sigma_X(X)}$. }

\dkz By restricting $(\ast)$ to the $x$th and $y$th coordinate
respectively, we obtain $k_{\sigma_X(X)}=l+l_{\sigma_X(X)}$ and
$k+k_{\sigma_X(X)}=l_{\sigma_X(X)}$, for $l=\sum_{b\in
B-\{\sigma_X(X)\}}l_b\cdot(\kappa_X(b))_x$, where
$(\kappa_X(b))_x$ is the $x$th coordinate of $\kappa_X(b)$, and
analogously, $k=\sum_{a\in
A-\{\sigma_X(X)\}}k_a\cdot(\kappa_X(a))_y$.

From this we infer $l+k=0$. Since $l,k\geq 0$, we obtain $l=k=0$,
which implies that $k_{\sigma_X(X)}=l_{\sigma_X(X)}$ and
$(\kappa_X(b))_x=(\kappa_X(a))_y=0$, for $b\in B-\{\sigma_X(X)\}$
and $a\in A-\{\sigma_X(X)\}$. So, $x$ and $y$ does not occur in
the elements of $B-\{\sigma_X(X)\}$ and $A-\{\sigma_X(X)\}$. \qed

\oprop{Lemma 6.3}{If $C$ is a simplicial complex such that
$X=\bigcup C$, and $B$ is a building set of $C$, then
$\kappa_X\!:TX\rightarrow \textbf{R}^{|X|}$ faithfully realizes
$C$ and $\St_X(C,B)$. Moreover, we have
$\bigcup\kappa_X(\St_X(C,B))\subseteq \bigcup\kappa_X (C)$. }

\dkz It is obvious that $\kappa_X$ faithfully realizes $C$ since
the set $\{\kappa_X(a)\mid a\in X\}$ is linearly independent. For
the rest of the proof we may rely on comments given in \cite{FY03}
(Section~5) immediately after the definition of $\Sigma({\cal
L},{\cal G})$. Here is how we prove these results.

Suppose $(\ast)$
$\sum_{a\in\alpha}k_a\cdot\kappa_X(a)=\sum_{b\in\beta}l_b\cdot
\kappa_X(b)$ for $\alpha,\beta\in\St_X(C,B)$ and $k_a,l_b>0$. Let
$\alpha',\beta'\in\Nt(C,B)$ be such that
$\alpha=P\sigma_X(\alpha')$ and $\beta=P\sigma_X(\beta')$. From
$(\ast)$ it follows that $x\in X$ occurs in some $a\in\alpha$ iff
it occurs in some $b\in\beta$. Hence, by Proposition 3.2 we
conclude that for
\[
\gamma=\{x\in X\mid x\;\mbox{\rm occurs in some }
a\in\alpha\}=\{x\in X\mid x\;\mbox{\rm occurs in some }
b\in\beta\}
\]
we have that $\alpha',\beta'\in\Nt(P(\gamma),B_\gamma)$, and hence
$\alpha,\beta\in\St_X(P(\gamma),B_\gamma)$.

To prove that $(\ast)$ implies that $\alpha=\beta$ and that
$k_a=l_a$ for every $a\in\alpha$, we proceed by induction on the
cardinality $|\gamma|$ of $\gamma$ by relying on the inductive
definition of a construction given in Section 2.3. (By Proposition
3.3 we have that $\alpha'$ and $\beta'$ are subsets of some
constructions of $B_\gamma$.) If $|\gamma|=0$, then we are done.

Suppose $|\gamma|\geq 1$ and $B_\gamma$ is connected, hence
$\gamma\in B_\gamma$. By the choice of $\gamma$, we have
$\bigcup\alpha'=\bigcup\beta'=\gamma$, and since every
$\alpha'$-antichain and every $\beta'$-antichain misses
$B_\gamma$, we must have that $\gamma\in\alpha'\cap\beta'$, and
hence $\sigma_X(\gamma)\in\alpha\cap\beta$. By the inductive
definition of a construction we have that
$\alpha-\{\sigma_X(\gamma)\}\in\St_X(P(\gamma-\{x\}),B_{\gamma-\{x\}})$
and
$\beta-\{\sigma_X(\gamma)\}\in\St_X(P(\gamma-\{y\}),B_{\gamma-\{y\}})$
for some $x,y\in\gamma$.

From Lemma 6.2 and Proposition 3.2, it follows that
\[
\alpha-\{\sigma_X(\gamma)\},\beta-\{\sigma_X(\gamma)\}\in
\St_X(P(\gamma-(\{x\}\cup\{y\})),B_{\gamma-(\{x\}\cup\{y\})})
\]
and that $k_{\sigma_X(\gamma)}=l_{\sigma_X(\gamma)}$. Hence, we
may cancel $k_{\sigma_X(\gamma)}\cdot\kappa_X(\sigma_X(\gamma))$
on both sides of $(\ast)$ and apply the induction hypothesis to
$\alpha-\{\sigma_X(\gamma)\}$, $\beta-\{\sigma_X(\gamma)\}$ and
$\gamma-(\{x\}\cup\{y\})$. (Note that
$|\gamma-(\{x\}\cup\{y\})|\leq|\gamma|-1$.)

Suppose $|\gamma|\geq 2$ and $B_\gamma$ is not connected. Then for
$\{B_{\gamma_1},\ldots,B_{\gamma_n}\}$, $n\geq 2$, being the
finest hypergraph partition of $B_\gamma$, we have that
$\alpha=\alpha_1\cup\ldots\cup\alpha_n$ and
$\beta=\beta_1\cup\ldots\cup\beta_n$, where
$\emptyset\neq\alpha_i,\beta_i\subseteq B_{\gamma_i}$. Since
$\{\kappa_X(x)\mid x\in\gamma\}$ is linearly independent, we have
that $(\ast)$ breaks into $n$ equations and we may apply the
induction hypothesis to each $\alpha_i$, $\beta_i$ and $\gamma_i$.
This concludes the first part of the proof.

To prove that  $\bigcup\kappa_X(\St_X(C,B))\subseteq
\bigcup\kappa_X (C)$, by Proposition 3.2, it is sufficient to show
that for every $\gamma\in C$, we have
\[
\bigcup \kappa_X(\St_X(P(\gamma),B_\gamma))\subseteq \mbox{\rm
cone}(\{\kappa_X(x)\mid x\in\gamma\})
\]
which is trivial. Actually, we may prove that
\[
\bigcup \kappa_X(\St_X(P(\gamma),B_\gamma))= \mbox{\rm
cone}(\{\kappa_X(x)\mid x\in\gamma\})
\]
relying on the inductive definition of a construction. This would
deliver that $\bigcup\kappa_X(\St_X(C,B))= \bigcup\kappa_X (C)$
but we don't need this stronger result here. \qed

\oprop{Remark 6.4}{If $\psi$ underlies an isomorphism between $C$
and $D$ and $f\cirk\psi$ faithfully realizes $C$,  then $f$
faithfully realizes $D$.}

Let $X$ be a finite, nonempty set and let $C$ be a simplicial
complex such that $\bigcup C\subseteq TX$. Let $B$ be a building
set of $C$. We define
\[
\Sb_X(C,B)=_{df}\{P\mu_X(\alpha)\mid \alpha\in\St_{TX}(C,B)\}.
\]
(It is easy to see that if $\bigcup C\subseteq X$, then
$\Sb_X(C,B)=\St_X(C,B)=\Ss_X(C,B)$.) The following proposition
shows when $\mu_X$ underlies an isomorphism between
$\St_{TX}(C,B)$ and $\Sb_X(C,B)$.

\prop{Proposition 6.5}{If $\kappa_X$ faithfully realizes the
simplicial complex $C$, then $\mu_X$ restricted to
$\bigcup\St_{TX}(C,B)$ is one-one and $\kappa_X$ faithfully
realizes $\Sb_X(C,B)$.}

\dkz Let $Y=\bigcup C$ and assume that $Y$ is linearly ordered in
an arbitrary way. By Lemma 6.3 we have that $\kappa_Y$ faithfully
realizes $C$ and $\St_Y(C,B)$. To simplify the notation, we assume
that $+^2$ is the formal addition of $TY$. Otherwise, we would
carry an extra isomorphism converting $+^2$ into the formal
addition of $TY$ at some places of this proof. By this assumption
we have that $\St_Y(C,B)=\St_{TX}(C,B)$, and hence $\kappa_Y$
faithfully realizes $\St_{TX}(C,B)$.

Let $L\!: \textbf{R}^{|Y|}\rightarrow \textbf{R}^{|X|}$ be a
linear transformation defined by $L\kappa_Y(a)=\kappa_X(a)$ for
$a\in Y$. Since $\kappa_X$ faithfully realizes $C$, we have that
$L\kappa_Y$ faithfully realizes $C$. By Lemma 6.3, it follows that
$\bigcup \kappa_Y(\St_{TX}(C,B))\subseteq \bigcup\kappa_Y(C)$.
Hence, the assumptions of Proposition 6.1 are satisfied for
$f=g=\kappa_Y$, $C$ and $D=\St_{TX}(C,B)$. So, we conclude that
$L\kappa_Y$ faithfully realizes $\St_{TX}(C,B)$.

By the definition of $L$ it follows that
$L\kappa_Y(c)=\kappa_X(\mu_X(c))$, for every $c\in TY$, and since
$\bigcup \St_{TX}(C,B)\subseteq TY$, this holds for every
$c\in\bigcup \St_{TX}(C,B)$. Hence, for $c,d\in\bigcup
\St_{TX}(C,B)$, we have that $\mu_X(c)=\mu_X(d)$ implies
$L\kappa_Y(c)=L\kappa_Y(d)$, which implies $c=d$ since $L\kappa_Y$
faithfully realizes $\St_{TX}(C,B)$. We conclude that $\mu_X$
restricted to $\bigcup\St_{TX}(C,B)$ is one-one.

By the above we have that $L\kappa_Y=\kappa_X\cirk\mu_X$
(restricted to $TY$), that $\mu_X$ underlies an isomorphism
between $\St_{TX}(C,B)$ and $\Sb_X(C,B)$, and that $L\kappa_Y$
faithfully realizes $\St_{TX}(C,B)$. Hence, by Remark 6.4,
$\kappa_X$ faithfully realizes $\Sb_X(C,B)$. \qed

The following remark is a corollary of Remark 5.4.

\prop{Remark 6.6}{If $\psi\!:X\rightarrow Y$ is a function that
underlies an isomorphism between simplicial complexes $C$ and $D$,
and $B$ is a building set of $C$, then the function $T\psi$
underlies an isomorphism between $\St_X(C,B)$ and
$\St_Y(D,P^2\psi(B))$.}

For $C$ a simplicial complex such that $\bigcup C\subseteq X$, and
$B_0,\ldots,B_n$ as in the definition of $\St_X(C,B_0,\ldots,B_n)$
(cf.\ Section~5), we have the following.

\prop{Proposition 6.7}{The restriction of $\mu_X^n$ to
$\bigcup\St_X(C,B_0,\ldots,B_n)$ is one-one and $\kappa_X$
faithfully realizes $\Ss_X(C,B_0,\ldots,B_n)$.}

\dkz We proceed by induction on $n$. In the basis, when $n=0$, we
have that $\mu_X^0$ is the identity function on $TX$, and by Lemma
6.3, $\kappa_X$ faithfully realizes $\St_X(C,B_0)$ which is equal
to $\Ss_X(C,B_0)$. For the induction step we use that
\[\St_X(C,B_0,\ldots,B_n)=\St_{T^nX}(\St_X(C,B_0,\ldots,B_{n-1}),P^2\sigma_X^n(B_n))\]
(cf.\ the end of the preceding section). By the induction
hypothesis and Remark 2.0.1, we have that
$\mu_X^{n-1}\!:T^nX\rightarrow TX$ underlies an isomorphism
between $\St_X(C,B_0,\ldots,B_{n-1})$ and
$C'=\Ss_X(C,B_0,\ldots,B_{n-1})$. Hence, by Remark 6.6,
$T\mu_X^{n-1}$ underlies an isomorphism between
$\St_X(C,B_0,\ldots,B_n)$ and $\St_{TX}(C',B)$, for
$B=P^2(\mu_X^{n-1}\cirk\sigma_X^n)(B_n)$. Since by the induction
hypothesis we have that $\kappa_X$ faithfully realizes $C'$, we
may apply Proposition 6.5 to conclude that $\mu_X$ restricted to
$\bigcup\St_{TX}(C',B)$ is one-one. So, the restriction of
$\mu_X^n=\mu_X\cirk T\mu_X^{n-1}$ to
$\bigcup\St_X(C,B_0,\ldots,B_n)$ is one-one.

By Proposition 6.5 we have that $\kappa_X$ faithfully realizes the
simplicial complex $\Sb_X(C',B)$. By Remark 2.0.1 we have that
$\mu_X^n$ underlies an isomorphism between
$\St_X(C,B_0,\ldots,B_n)$ and $\Sb_X(C',B)$, and between
$\St_X(C,B_0,\ldots,B_n)$ and $\Ss_X(C,B_0,\ldots,B_n)$. Hence, by
Remark 2.0.2, we have that $\Sb_X(C',B)=\Ss_X(C,B_0,\ldots,B_n)$.
So, we obtain that $\kappa_X$ faithfully realizes the simplicial
complex $\Ss_X(C,B_0,\ldots,B_n)$. \qed

As an immediate corollary of Propositions 5.3 and 6.7 we have the
following for $C$, $X$ and $B_0,\ldots,B_n$ as above.

\propo{Theorem 6.8}{The function $\mu_X^n\cirk\sigma_X^{n+1}$
underlies an isomorphism between $\Nt(C,B_0,\ldots,B_n)$ and
$\Ss_X(C,B_0,\ldots,B_n)$. }

\section{An alternative approach}

Although Theorem 6.8 provides us with an efficient notation for
the result of the iteration of $\Nt$, still we have to calculate
first $\Nt(C,B_0,\ldots,B_n)$ in order to obtain
$\Ss_X(C,B_0,\ldots,B_n)$. Also, it is not easy to write down
$B_n$ for large $n$, and since we always have that
$\{\{\beta\}\mid \beta\in B_i\}\subseteq B_{i+1}$, there is some
superfluous information carried in this notation. The aim of this
section is to formulate a direct procedure of calculating these
simplicial complexes and to give a simple, more economical
notation for building sets. We will just formulate here the
results without going into the proofs, which are straightforward
but tedious.

Let $X$ be a finite and nonempty set. For $\beta\in P_F(TX)$, let
$\beta^+=_{df}\mu_X\sigma_{TX}(\beta)$. For example, if
$\beta=\{x,x+y\}$, then $\beta^+=2x+y$. For $B\subseteq P_F(TX)$,
let $[B]^+=P(\mu_X\sigma_{TX})(B)=\{\beta^+\mid\beta\in B\}$.

Let $\alpha\subseteq TX$ be a finite set such that $\kappa_X$
faithfully realizes the simplicial complex $P(\alpha)$. This means
that $\{\kappa_X(a)\mid a\in\alpha\}$ is linearly independent.
Hence, $\mu_X\sigma_{TX}$ is one-one on $P_F(\alpha)$ and
$P(\mu_X\sigma_{TX})$ is one-one on $P(P_F(\alpha))$.

We say that
$D\subseteq[P_F(\alpha)]^+=[P(\alpha)-\{\emptyset\}]^+$ is a
\emph{flat building set} of $P(\alpha)$ when the following
conditions hold:
\begin{tabbing}
\hspace{1.5em}\=(D1)\hspace{1em}\=If
$\gamma,\delta\subseteq\alpha$ and $\gamma\cap \delta\neq
\emptyset$ and $\gamma^+,\delta^+\in D$, then $(\gamma\cup
\delta)^+\in D$;
\\[1ex]
\>(D2)\>$D\cap \alpha=\emptyset$.
\end{tabbing}
This means that there is a building set $B$ of $P(\alpha)$ such
that $[B]^+ -\alpha=D$. For example, if $\alpha=\{x,x+y,x+y+z\}$,
then $\{2x+y\}$ is a flat building set of $P(\alpha)$.

For $D$ a flat building set of $P(\alpha)$ and
$\beta\subseteq\alpha$, let $D_\beta=D\cap[P_F(\beta)]^+$ (cf.\
Section 2.3). If $B$ is a building set of $P(\alpha)$ such that
$[B]^+ -\alpha=D$, then, since $\mu_X\sigma_{TX}$ is one-one, we
have that $D_\beta=[B_\beta]^+-\beta$. So, $D_\beta$ is a flat
building set of $P(\beta)$. By induction on the cardinality $n\geq
0$ of $D$, we define the set $\C_X(\alpha,D)\subseteq P_F(TX)$ of
\emph{flat constructions} of a flat building set $D$ of
$P(\alpha)$ as follows

\begin{itemize}
\item[(0)] if $n=0$, i.e.\ $D=\emptyset$, then
$\C_X(\alpha,D)=\{\alpha\}$;

\item[(1)] if $n\geq 1$ and $\alpha^+\in D$, then
\[
\C_X(\alpha,D)=\{\{\alpha^+\}\cup\gamma\mid (\exists
x\in\alpha)\;\gamma\in\C_X(\alpha-\{x\},D_{\alpha-\{x\}})\};
\]

\item[(2)] if $n\geq 1$ and $\alpha^+\not\in D$, then, for
$M=\max\{\beta\mid\beta\subseteq\alpha\; \&\; \beta^+\in
D\cup\alpha\}$,
\[
\C_X(\alpha,D)=\{\gamma _1\cup\ldots\cup\gamma_n\mid
(\gamma_1,\ldots,\gamma_n)\in\prod_{\beta\in M}
\C_X(\beta,D_\beta)\}.
\]
\end{itemize}
If $\alpha=\{x,x+y,x+y+z\}$ and $D=\{2x+y\}$, then
$\C_X(\alpha,D)=\{\{x,2x+y,x+y+z\},\{x+y,2x+y,x+y+z\}\}$.

For a building set $B$ of $P(\alpha)$, by comparing the definition
of a construction given in Section 2.3 and the definition of a
flat construction given above we have the following.

\prop{Proposition 7.1}{$K$ is a construction of $B$ iff $[K]^+$ is
a flat construction of $[B]^+ -\alpha$.}

Let now $C$ be a simplicial complex such that $\bigcup C\subseteq
TX$ and suppose that $\kappa_X$ faithfully realizes $C$. We say
that $D\subseteq \bigcup\{[P_F(\alpha)]^+\mid\alpha\in C\}$ is a
\emph{flat building set} of $C$, when for every $\alpha\in C$ we
have that $D_\alpha=D\cap[P_F(\alpha)]^+$ is a flat building set
of $P(\alpha)$. This means that there is a building set $B$ of $C$
such that $[B]^+ -\bigcup C=D$. Let $\Et_X(C,D)$ be the simplicial
complex based on
\[
\bigcup\{\C_X(\alpha,D_\alpha)\mid\alpha\; \mbox{\rm is a basis of
}C\}.
\]
For $C$ and $D$ as above, let $B$ be a building set of $C$ such
that $[B]^+ -\bigcup C=D$. From the definition given after Remark
6.4, we conclude that $\Sb_X(C,B)=\{[N]^+\mid N\in\Nt(C,B)\}$, and
by relying on Proposition 7.1, we have the following.

\prop{Proposition 7.2}{$\Sb_X(C,B)=\Et_X(C,D)$.}

From Proposition 6.5 we conclude that $\kappa_X$ faithfully
realizes $\Et_X(C,D)$. Since $\bigcup\Et_X(C,D)\subseteq TX$, we
may iterate this procedure and for $n\geq 0$, $D_0$ a flat
building set of $C$, and $D_{i+1}$, for $0\leq i\leq n$ a flat
building set of $\Et_X(C,D_0,\ldots,D_i)$, we define
\[
\Et_X(C,D_0,\ldots,D_{n+1})=_{df}\Et_X(\Et_X(C,D_0,\ldots,D_n),D_{n+1}).
\]

Let $C$ be a simplicial complex such that $\bigcup C\subseteq X$,
and let $B_0,\ldots,B_n$ be as in the definition of
$\St(C,B_0,\ldots,B_n)$. Let $G_0=B_0$, and for $1\leq i\leq n$,
let $G_i=P^2(\mu_X^{i-1}\cirk\sigma_X^i)(B_i)$. Let $D_0=[G_0]^+
-\bigcup C$, and for $1\leq i\leq n$, let $D_i=[G_i]^+ -
[G_{i-1}]^+$. By using Proposition 7.2 and results of the
preceding section, we have the following.

\prop{Proposition 7.3}{$\Ss_X(C,B_0,\ldots,B_n)
=\Et_X(C,D_0,\ldots,D_n)$.}

In the table below we give, as an example, the bases of simplicial
complexes $C$, $\Et_X(C,D_0)$, $\Et_X(C,D_0, D_1)$ and
$\Et_X(C,D_0,D_1,D_2)$, where \linebreak $X=\{x,y,z,u\}$,
$C=P(X)-\{X\}$, $D_0=\{x+y,x+y+z\}$,
$D_1=\{2x+y,2x+y+z,2x+2y+z,3x+2y+z\}$ and
$D_2=\{6x+5y+3z,3x+3y+2z\}$.

\begin{table}
\begin{center}
\begin{tabular}{|l|l|l|l|} \hline & & & \\[-2.2ex]
$C$ & $\Et_X(C,D_0)$ & $\Et_X(C,D_0, D_1)$ &
$\Et_X(C,D_0,D_1,D_2)$
\\[.3ex] \hline\hline & & & \\[-2.2ex] $\{x,y,z\}$ & $\{x,x+y,$ &
$\{x+y+z,2x+2y+z,$ & $\{x+y+z,3x+3y+2z,$ \\
& $x+y+z\}$ & $3x+2y+z\}$ & $6x+5y+3z\}$
\\[.3ex] \cline{4-4} & & & \\[-2.2ex] & & & $\{2x+2y+z,3x\pl 3y\pl 2z,$ \\ & & & $6x+5y+3z\}$
\\[.3ex] \cline{4-4} & & & \\[-2.2ex] & & & $\{x+y+z,3x+2y+z,$ \\ & & & $6x+5y+3z\}$
\\[.3ex] \cline{4-4} & & & \\[-2.2ex] & & & $\{2x+2y+z,3x+2y+z,$ \\ & & & $6x+5y+3z\}$
\\[.3ex] \cline{3-4} & & & \\[-2.2ex] & & $\{x+y+z,2x+y+z,$ & $\{x+y+z,2x+y+z,$ \\ & & $3x+2y+z\}$ & $3x+2y+z\}$
\\[.3ex] \cline{3-4} & & & \\[-2.2ex] & & $\{x+y,2x+2y+z,$ & $\{x+y,2x+2y+z,$ \\ & & $3x+2y+z\}$ & $3x+2y+z\}$
\\[.3ex] \cline{3-4} & & & \\[-2.2ex] & & $\{x+y,2x+y,$ & $\{x+y,2x+y,$ \\ & & $3x+2y+z\}$ & $3x+2y+z\}$
\\[.3ex] \cline{3-4} & & & \\[-2.2ex] & & $\{x,2x+y+z,$ & $\{x,2x+y+z,$ \\ & & $3x+2y+z\}$ & $3x+2y+z\}$
\\[.3ex] \cline{3-4} & & & \\[-2.2ex] & & $\{x,2x+y,3x+2y+z\}$ & $\{x,2x+y,3x+2y+z\}$
\\[.3ex] \cline{2-4} & & & \\[-2.2ex] & $\{y,x+y,$ & $\{y,x+y+z,$ & $\{y,x+y+z,$ \\ & $x+y+z\}$ & $2x+2y+z\}$ & $3x+3y+2z\}$
\\[.3ex] \cline{4-4} & & & \\[-2.2ex] & & & $\{y,2x+2y+z,$ \\ & & & $3x+3y+2z\}$
\\[.3ex] \cline{3-4} & & & \\[-2.2ex] & & $\{y,x+y,2x+2y+z\}$ & $\{y,x+y,2x+2y+z\}$
\\[.3ex] \cline{2-4} & & & \\[-2.2ex] & $\{x,z,$ & $\{x,z,2x+y+z\}$ & $\{x,z,2x+y+z\}$
\\ \cline{3-4} & $x+y+z\}$ & & \\[-1.7ex] & & $\{z,x+y+z,2x+y+z\}$ & $\{z,x+y+z,2x+y+z\}$
\\[.3ex] \cline{2-4} & & & \\[-2.2ex] & $\{y,z,$ & $\{y,z,x+y+z\}$ & $\{y,z,x+y+z\}$
\\ & $x+y+z\}$ & &
\\[.3ex] \hline & & & \\[-2.2ex] $\{x,y,u\}$ & $\{x,u,x\pl y\}$ &
$\{x,u,2x+y\}$ & $\{x,u,2x+y\}$
\\[.3ex] \cline{3-4} & & & \\[-2.2ex] & & $\{u,x+y,2x+y\}$ & $\{u,x+y,2x+y\}$
\\[.3ex] \cline{2-4} & & & \\[-2.2ex] & $\{y,u,x\pl y\}$ & $\{y,u,x+y\}$ & $\{y,u,x+y\}$
\\[.3ex] \hline & & & \\[-2.2ex] $\{x,z,u\}$ & $\{x,z,u\}$ &
$\{x,z,u\}$ & $\{x,z,u\}$
\\[.3ex] \hline & & & \\[-2.2ex] $\{y,z,u\}$ & $\{y,z,u\}$ &
$\{y,z,u\}$ & $\{y,z,u\}$
\\[.3ex] \hline

\end{tabular}
\end{center}
\end{table}

\section{$\Sb_X$ and combinatorial blowups}

Feichtner and Kozlov defined in \cite{FK04} (Definition 3.1; see
also \cite{C07}, Definition 1.4) the poset $\Bl_\alpha{\cal L}$,
called the \emph{combinatorial blowup} of a finite-meet
semilattice $\cal L$ at its element $\alpha$, in the following
way. (Here the authors presumably assumed that $\alpha$ is not the
bottom element of $\cal L$.) The set of elements of
$\Bl_\alpha{\cal L}$ is
\[
\{\gamma\in{\cal L}\mid
\gamma\not\geq\alpha\}\cup\{(\alpha,\gamma)\mid \gamma\in{\cal
L}\;\mbox{\rm and } \gamma\not\geq\alpha\;\mbox{\rm and
}(\gamma\vee\alpha)_{\cal L}\;\mbox{\rm exists}\},
\]
while for $\beta,\gamma\not\geq\alpha$, the order relation is
given by
\begin{tabbing}
\hspace{1.5em}\=(1)\quad\=$\beta>\gamma$ in $\Bl_\alpha{\cal L}$
if $\beta>\gamma$ in $\cal L$;
\\[1ex]
\>(2)\>$(\alpha,\beta)>(\alpha,\gamma)$ in $\Bl_\alpha{\cal L}$ if
$\beta>\gamma$ in $\cal L$;
\\[1ex]
\>(3)\>$(\alpha,\beta)>\gamma$ in $\Bl_\alpha{\cal L}$ if
$\beta\geq\gamma$ in $\cal L$.
\end{tabbing}

Let $X$ be finite, nonempty set, and let $C$ be a simplicial
complex such that $Y=\bigcup C\subseteq TX$ and $\kappa_X$
faithfully realizes $C$. Note that for every $\alpha\in C$ that is
neither empty nor a singleton, we have that
$\alpha^+=\mu_X\sigma_{TX}(\alpha)\not\in Y$. It is easy to see
that for such $C$, we have that $\Bl_\alpha C$, for
$\alpha\neq\emptyset$, is isomorphic to the simplicial complex
\[
\{\gamma\in
C\mid\alpha\not\subseteq\gamma\}\cup\{\gamma\cup\{\alpha^+\}\mid\alpha\not\subseteq\gamma\;
\mbox{\rm and }\alpha\cup\gamma\in C\},
\]
which we will also denote by $\Bl_\alpha C$ and call the
\emph{blowup} of $C$ at $\alpha$. Note that if $\alpha$ is a
singleton, then $\Bl_\alpha C=C$.

For $X$ and $C$ as above we prove the following proposition whose
content may be derived from the proof of Proposition 2.4 of
\cite{C07}.

\prop{Proposition 8.1}{If $B$ and $B\cup\{\beta\}$ are building
sets of $C$, then for $\alpha=[\max B_\beta]^+$ we have
\[
\Sb_X(C,B\cup\{\beta\})=\Bl_\alpha\, \Sb_X(C,B).
\]
}

\dkz If the cardinality $|\max B_\beta|$ of $\max B_\beta$ is 1,
then  since $B$ is a building set of $C$, we have that $\beta\in
B$. Hence, $\Sb_X(C,B\cup\{\beta\})=\Sb_X(C,B)=\Bl_\alpha\,
\Sb(C,B)$, since $\alpha$ is a singleton. So, we may proceed with
the assumption that $|\max B_\beta|\geq 2$, i.e.\ $\beta\not\in
B$.

For the proof of $\subseteq$-direction, let $\gamma \in
\Sb_X(C,B\cup\{\beta\})$. Note that
$\Sb_X(C,B\cup\{\beta\})=\{[N]^+\mid N\in\Nt(C,B\cup\{\beta\})\}$.
According to Proposition 3.2, we have that $\gamma=[N]^+$ for some
$N$ that satisfies $N\subseteq B\cup\{\beta\}$, $\bigcup N\in C$
and every $N$-antichain misses $B\cup\{\beta\}$. Since $B$ is a
building set of $C$, we have that the members of $\max B_\beta$
are mutually disjoint and $\bigcup(\max B_\beta)=\beta$. We may
conclude that $\alpha^+=\beta^+$ and $\alpha\not\subseteq \gamma$
(otherwise, since $\kappa_X$ faithfully realizes $C$, $\max
B_\beta$ would be an $N$-antichain that does not miss
$B\cup\{\beta\}$).

Suppose $\beta^+\not\in\gamma$, which is equivalent to
$\beta\not\in N$ since $\kappa_X$ faithfully realizes $C$. We have
that $N\subseteq B$ and every $N$-antichain misses $B$ since it
misses $B\cup\{\beta\}$. By Proposition 3.2, we have that
$\gamma\in\Sb_X(C,B)$, and since $\alpha\not\subseteq\gamma$, we
conclude that $\gamma\in\Bl_\alpha\,\Sb_X(C,B)$.

Suppose $\beta^+\in\gamma$. We want to show that
$\alpha\cup(\gamma-\{\alpha^+\})\in\Sb_X(C,B)$. Since
$\alpha^+=\beta^+$, we have that
$\gamma-\{\alpha^+\}=\gamma-\{\beta^+\}$, and hence
$\alpha\cup(\gamma-\{\alpha^+\})=[M]^+$ for $M=\max
B_\beta\cup(N-\{\beta\})$. We have that $M\subseteq B$ and
$\bigcup M=\bigcup N\in C$. We have to show that every
$M$-antichain misses $B$.

Suppose $S$ is an $M$-antichain that does not miss $B$. We have
that $S\cap \max B_\beta\neq\emptyset$, since otherwise $S$ would
be an $N$-antichain that does not miss $B\cup\{\beta\}$. Let
$\beta'\in S\cap\max B_\beta$ and let $S'=S-B_\beta$.

If $S'=\emptyset$, then $\bigcup S\subseteq \beta$. From $\bigcup
S\in B$ we conclude that $\beta'\subseteq\bigcup S\in B_\beta$ and
hence $\beta'=\bigcup S$, which together with $\beta'\in S$
contradicts the assumption that $S$ is an $M$-antichain.

If $S'\neq\emptyset$, then for every $\gamma\in S'$ we have
$\beta\not\subseteq\gamma$. Otherwise, $\beta',\gamma\in S$ and
$\beta'\subseteq\gamma$ which contradicts the assumption that $S$
is an $M$-antichain. Since $\gamma\not\subseteq\beta$ holds by the
definition of $S'$, we have that $S'\cup\{\beta\}$ is an
$N$-antichain. We have that $(\bigcup S')\cup\beta=(\bigcup
S)\cup\beta\subseteq\bigcup N\in C$, and from
$\emptyset\neq\beta'\subseteq(\bigcup S)\cap\beta$, $\bigcup S\in
B\cup\{\beta\}$ and $\beta\in B\cup\{\beta\}$, since
$B\cup\{\beta\}$ is a building set of $C$, we have that $(\bigcup
S)\cup\beta\in B\cup\{\beta\}$. So, $S'\cup\{\beta\}$ is an
$N$-antichain that does not miss $B\cup\{\beta\}$, which is a
contradiction.

By Proposition 3.2 we have that
$\alpha\cup(\gamma-\{\alpha^+\})\in\Sb_X(C,B)$. Since
$\alpha\not\subseteq\gamma$, we have that
$\alpha\not\subseteq\gamma-\{\alpha^+\}$, and hence
$(\gamma-\{\alpha^+\})\cup\{\alpha^+\}=\gamma\in\Bl_\alpha\,\Sb_X(C,B)$.

For the proof of $\supseteq$-direction, suppose first that
$\gamma\in\Bl_\alpha\, \Sb_X(C,B)$ is such that
$\gamma\in\Sb_X(C,B)$ and $\alpha\not\subseteq\gamma$. By
Proposition 3.2, we have that $\gamma=[N]^+$ for some $N$ that
satisfies $N\subseteq B$, $\bigcup N\in C$ and every $N$-antichain
misses $B$.

We show that
\begin{tabbing}
\hspace{1.5em}$(\ast)$\quad if $S$ is an $N$-antichain, then
$\bigcup S\neq\beta$.
\end{tabbing}
Since $\alpha\not\subseteq\gamma$, we have that $\max
B_\beta\not\subseteq N$. So, there is $\beta'\in\max B_\beta$ such
that $\beta'\not\in N$. Suppose $\bigcup S=\beta$, then since all
the elements of $\max B_\beta$ are mutually disjoint and
$S\subseteq B_\beta$, there must be a subset $S'$ of $S$ such that
$\bigcup S'=\beta'$. From $\beta'\not\in N$ we conclude that
$\beta'\not\in S'$ and hence $|S'|\geq 2$. So, $S'$ is an
$N$-antichain that does not miss $B$, which is a contradiction.
Hence, $(\ast)$ holds and this guarantees that every $N$-antichain
misses $B\cup\{\beta\}$. By Proposition 3.2 we have that
$\gamma\in\Sb_X(C,B\cup\{\beta\})$.

Suppose now that $\gamma\cup\{\alpha^+\}\in\Bl_\alpha\,
\Sb_X(C,B)$, where $\alpha\cup\gamma\in\Sb_X(C,B)$ and
$\alpha\not\subseteq\gamma$. Since $\Sb_X(C,B)$ is a simplicial
complex, we have $\gamma\in\Sb_X(C,B)$. By Proposition 3.2, we
have that $\gamma=[N]^+$ and $\alpha\cup\gamma=[M]^+$ for $N$ and
$M$ such that $N,M\subseteq B$, $\bigcup N,\bigcup M\in C$ and
every $N$-antichain and every $M$-antichain misses $B$.

Let $K=N\cup\{\beta\}$. We have that $\gamma\cup\{\alpha^+\}=
[K]^+$ and $K\subseteq B\cup\{\beta\}$ and $\bigcup K=(\bigcup
N)\cup\beta=(\bigcup N)\cup \max B_\beta=\bigcup M\in C$. We have
to show that every $K$-antichain misses $B\cup\{\beta\}$.

If $S$ is a $K$-antichain such that $\beta\not\in S$, then $S$ is
an $N$-antichain and it misses $B$. We can prove that $\bigcup
S\neq\beta$ as we proved $(\ast)$ above, and hence $S$ is a
$K$-antichain that misses $B\cup\{\beta\}$.

If $S$ is a $K$-antichain such that $\beta\in S$, then $\bigcup
S\neq\beta$ (otherwise, $S$ is not a $K$-antichain). If $\bigcup
S\in B$, then for $S'=(S-\{\beta\})\cup \max B_\beta$ we have
$\bigcup S'=\bigcup S\in B$. By eliminating from $S'$ every member
properly contained in some other member of $S'$, we obtain an
$M$-antichain $S''$ such that $\bigcup S''\in B$ (it is easy to
see that $|S''|\geq 2$), which is a contradiction. So $S$ is a
$K$-antichain that misses $B\cup\{\beta\}$.

Hence, every $K$-antichain misses $B\cup\{\beta\}$, and by
Proposition 3.2, we have that
$\gamma\cup\{\alpha^+\}\in\Sb_X(C,B\cup\{\beta\})$. \qed

Let $X$ be finite, nonempty set, and let $C$ be a simplicial
complex such that $\emptyset\neq Y=\bigcup C\subseteq TX$ and
$\kappa_X$ faithfully realizes $C$. Let $B$ be a building set of
$C$. Then we have the following.

\prop{Proposition 8.2}{There exist $\beta_1,\ldots,\beta_m\in B$
such that $\Sb_X(C,B)=\Bl_{\beta_1}(\ldots\Bl_{\beta_m}\, C)$}

\dkz We proceed by induction on the cardinality of $B-B_\bot$,
where $B_\bot=\{\{a\}\mid a\in Y\}$. If $B=B_\bot$, then
$\Sb_X(C,B)=C=\Bl_{\{a\}}\, C$ for arbitrary $a\in Y$. Otherwise,
let $\beta_1\in\min (B-B_\bot)$. By (5) of Proposition 3.1 we have
that $B-\{\beta_1\}$ is a building set of $C$. By Proposition 8.1,
since $[\max B_{\beta_1}]^+=\beta_1$, we have that
$\Sb_X(C,B)=\Bl_{\beta_1}\, \Sb_X(C,B-\{\beta_1\})$ and we may
apply the induction hypothesis. \qed

For $C$ a simplicial complex such that $\emptyset\neq\bigcup
C\subseteq X$, and $B_0,\ldots,B_n$ as in the definition of
$\St_X(C,B_0,\ldots,B_n)$ (cf.\ Section~5), we have the following.

\prop{Proposition 8.3}{There exist $\beta_1,\ldots,\beta_m\in
P_F(TX)$ such that \linebreak $\Ss_X(C,B_0,\ldots,
B_n)=\Bl_{\beta_1}(\ldots\Bl_{\beta_m}\, C)$. }

\dkz We proceed by induction on $n$. If $n=0$, then
$\Ss_X(C,B_0)=\Sb_X(C,B_0)$ and we may apply Proposition 8.2,
since $\kappa_X$ faithfully realizes $C$. If $n\geq 1$, then as in
the proof of Proposition 6.7, we have that
\[
\Ss_X(C,B_0,\ldots, B_n)=\Sb_X(\Ss_X(C,B_0,\ldots, B_{n-1}),B),
\]
for $B=P^2(\mu_X^{n-1}\cirk\sigma_X^n)(B_n)$, and $\kappa_X$
faithfully realizes $\Ss_X(C,B_0,\ldots, B_{n-1})$. By the
induction hypothesis $\Ss_X(C,B_0,\ldots,
B_{n-1})=\Bl_{\beta_1}(\ldots\Bl_{\beta_k}\, C)$ and it remains
just to apply Proposition 8.2 to
$\Sb_X(\Bl_{\beta_1}(\ldots\Bl_{\beta_k}\, C),B)$. \qed

\section{Stellar subdivision and truncation of polytopes}

In this section we elaborate two operations on polytopes, one dual
to the other, which realize the operation of combinatorial blowup
on the face lattices of polytopes. For a general reference to the
theory of polytopes, we refer the reader to \cite{G03} and
\cite{Z95}. We also try to keep to the notation used in these two
books.

According to \cite{HRZ04} (see also \cite{J07}) a stellar
subdivision of a polytope $P$ in a proper face $F$ is a polytope
$\conv\,(P\cup\{x^F\})$ where $x^F$ is a point of the form
$y^F-\varepsilon(y^P-y^F)$, where $y^P$ is in the interior of $P$,
$y^F$ is in the relative interior of $F$, and $\varepsilon$ is
small enough. We use $\st_F\, P$ to denote a stellar subdivision
of $P$ in $F$. For example, if $P$ is a cube $A_1B_1C_1D_1A_2 B_2
C_2 D_2$ given on the left-hand side, then $\st_{B_1B_2}\, P$ is
illustrated on the right-hand side:

\begin{center}
\begin{picture}(300,120)
{\thicklines \put(0,10){\line(1,0){70}} \put(0,90){\line(1,0){70}}
\put(30,110){\line(1,0){70}} \put(0,10){\line(0,1){80}}
\put(100,30){\line(0,1){80}} \put(0,90){\line(3,2){30}}
\put(70,90){\line(3,2){30}} \put(70,10){\line(3,2){30}}


\multiput(0,10)(3,2){10}{\makebox(0,0){\circle*{.5}}}
\multiput(30,30)(3,0){24}{\makebox(0,0){\circle*{.5}}}
\multiput(30,30)(0,3){27}{\makebox(0,0){\circle*{.5}}}

\put(70,10){\line(0,1){80}}}

\put(-1,7){\makebox(0,0)[tr]{$A_1$}}
\put(70,7){\makebox(0,0)[tr]{$B_1$}}
\put(102,30){\makebox(0,0)[bl]{$C_1$}}
\put(28,30){\makebox(0,0)[br]{$D_1$}}
\put(0,92){\makebox(0,0)[br]{$A_2$}}
\put(70,92){\makebox(0,0)[br]{$B_2$}}
\put(101,111){\makebox(0,0)[bl]{$C_2$}}
\put(28,111){\makebox(0,0)[br]{$D_2$}}


{\thicklines \put(200,10){\line(1,0){70}}
\put(200,90){\line(1,0){70}} \put(230,110){\line(1,0){70}}
\put(200,10){\line(0,1){80}} \put(300,30){\line(0,1){80}}
\put(200,90){\line(3,2){30}} \put(270,90){\line(3,2){30}}
\put(270,10){\line(3,2){30}}

\put(200,10){\line(2,1){80}} \put(200,90){\line(2,-1){80}}
\put(270,10){\line(1,4){10}} \put(270,90){\line(1,-4){10}}
\put(300,30){\line(-1,1){20}} \put(300,110){\line(-1,-3){20}}

\multiput(200,10)(3,2){10}{\makebox(0,0){\circle*{.5}}}
\multiput(230,30)(3,0){24}{\makebox(0,0){\circle*{.5}}}
\multiput(230,30)(0,3){27}{\makebox(0,0){\circle*{.5}}}}


\put(199,7){\makebox(0,0)[tr]{$A_1$}}
\put(270,7){\makebox(0,0)[tr]{$B_1$}}
\put(302,30){\makebox(0,0)[bl]{$C_1$}}
\put(228,30){\makebox(0,0)[br]{$D_1$}}
\put(200,92){\makebox(0,0)[br]{$A_2$}}
\put(270,92){\makebox(0,0)[br]{$B_2$}}
\put(301,111){\makebox(0,0)[bl]{$C_2$}}
\put(228,111){\makebox(0,0)[br]{$D_2$}}
\put(284,47){\makebox(0,0)[bl]{$B$}}

\end{picture}
\end{center}

To show that the stellar subdivision is well defined, i.e.\ that
the face lattice of $\st_F\, P$ does not depend on the choice of
$x^F$, we need a more precise definition of this notion. For this
we rely on some notions introduced in \cite{G03}.

Let $P\subseteq \R^d$ be a $d$-polytope, $H$ a hyperplane such
that $H\cap \int\, P=\emptyset$, and let $V\in\R^d$. Then $V$ is
\emph{beneath}, or \emph{beyond} $H$ (with respect to $P$),
provided $V$ belongs to the open halfspace determined by $H$ which
contains $\int\, P$, or does not meet $P$, respectively. If $V\in
\R^d$ and $F$ is a facet of the $d$-polytope $P\subseteq \R^d$,
then $V$ is beneath $F$ or beyond $F$ provided $V$ is beneath or
beyond $\aff\, F$, respectively. A \emph{stellar subdivision}
$\st_F\, P$ of a polytope $P$ in a proper face $F$ is a polytope
$\conv(P\cup\{x^F\})$ where $x^F$ is a point beneath every facet
not containing $F$ and beyond every facet containing $F$.

Let $\F(P)$ denote the face lattice of a polytope $P$ and let
$\F^-(P)$ denote the meet-semilattice obtained from $\F(P)$ by
removing the top element $P$. The following proposition connects
the operation of combinatorial blowup on finite-meet semilattices
and the operation of stellar subdivision on polytopes.

\prop{Proposition 9.1}{For a proper face $F$ of a polytope $P$ we
have that
\[
\F^-(\st_F\, P)\cong\Bl_F\, \F^-(P).
\]}

\dkz Apply Theorem~1 of \cite{G03} (Section 5.2).
\qed

We define now an operation that is polar (dual) to stellar
subdivision, which we will call truncation. This operation is
mentioned in \cite{G03} under the name ``cutting off'' and it
appears, under different names, in some recent publications
(``blow-up'' in \cite{CD06}, ``shaving construction'' in
\cite{PRW08}, etc.)

Let $P\subseteq \R^d$ be a $d$-polytope, $V\in P$, and let $\pi^+$
be a halfspace. We say that $\pi^+$ is \emph{beneath} $V$ when $V$
belongs to $\int\, \pi^+$, and we say that $\pi^+$ is
\emph{beyond} $V$ when $V$ does not belong to $\pi^+$.

A \emph{truncation} $\tr_F\, P$ of a polytope $P$ in a proper face
$F$ is a polytope $P\cap\pi^+$ where $\pi^+$ is a halfspace
beneath every vertex not contained in $F$ and beyond every vertex
contained in $F$. This defines an operation dual to stellar
subdivision and hence we have the following proposition.

\propo{Proposition 9.2}{For a proper face $F$ of a polytope $P$ we
have that (up to combinatorial equivalence)
\[
(\st_F\, P)^\triangle=\tr_{F^\diamond}\, P^\triangle.
\]}

\noindent For example, if $P$ is an octahedron $ABCDMN$ given on
the left-hand side, then $\tr_{AB}\, ABCDMN$ is illustrated on the
right-hand side:

\begin{center}
\begin{picture}(300,170)
{\thicklines \put(20,40){\line(1,-1){30}}
\put(80,40){\line(-1,-1){30}}

\put(10,80){\line(1,2){40}} \put(90,80){\line(-1,2){40}}
\put(0,70){\line(0,1){30}} \put(100,70){\line(0,1){30}}
\put(0,100){\line(5,6){50}} \put(100,100){\line(-5,6){50}}

\multiput(1,100)(3,0){34}{\makebox(0,0){\circle*{.5}}}
\multiput(0.5,100)(2.5,-4.5){20}{\makebox(0,0){\circle*{.5}}}
\multiput(99.5,100)(-2.5,-4.5){20}{\makebox(0,0){\circle*{.5}}} }

{\thicklines \put(0,70){\line(0,-1){10}}
\put(100,70){\line(0,-1){10}} \put(0,60){\line(1,0){100}}
\put(0,60){\line(1,2){10}} \put(100,60){\line(-1,2){10}}
\put(0,60){\line(1,-1){20}} \put(100,60){\line(-1,-1){20}} }

\put(-2,60){\makebox(0,0)[r]{$A$}}
\put(102,60){\makebox(0,0)[l]{$B$}}
\put(102,100){\makebox(0,0)[l]{$C$}}
\put(-2,100){\makebox(0,0)[r]{$D$}}
\put(50,163){\makebox(0,0)[b]{$M$}}
\put(50,7){\makebox(0,0)[t]{$N$}}


{\thicklines \put(250,10){\line(1,1){40}}
\put(250,10){\line(-1,1){40}} \put(210,50){\line(1,0){80}}
\put(210,80){\line(1,0){80}} \put(210,50){\line(-2,3){10}}
\put(290,50){\line(2,3){10}} \put(210,80){\line(-2,-3){10}}
\put(290,80){\line(2,-3){10}} \put(210,80){\line(1,2){40}}
\put(290,80){\line(-1,2){40}} \put(200,65){\line(0,1){35}}
\put(300,65){\line(0,1){35}} \put(200,100){\line(5,6){50}}
\put(300,100){\line(-5,6){50}}

\multiput(201,100)(3,0){34}{\makebox(0,0){\circle*{.5}}}
\multiput(200.5,100)(2.5,-4.5){20}{\makebox(0,0){\circle*{.5}}}
\multiput(299.5,100)(-2.5,-4.5){20}{\makebox(0,0){\circle*{.5}}} }

\put(210,47){\makebox(0,0)[r]{$A_1$}}
\put(290,47){\makebox(0,0)[l]{$B_1$}}
\put(200,62){\makebox(0,0)[r]{$A_2$}}
\put(300,62){\makebox(0,0)[l]{$B_2$}}
\put(230,87){\makebox(0,0)[r]{$A_3$}}
\put(270,87){\makebox(0,0)[l]{$B_3$}}
\put(302,100){\makebox(0,0)[l]{$C$}}
\put(198,100){\makebox(0,0)[r]{$D$}}
\put(250,163){\makebox(0,0)[b]{$M$}}
\put(250,7){\makebox(0,0)[t]{$N$}}

\end{picture}
\end{center}

Let $\F_-(P)$ denote the join-semilattice obtained from $\F(P)$ by
removing the bottom element $\emptyset$. Then we have the
following proposition, which is analogous to Proposition 9.1 (the
operation $^{op}$ reverses the order of a poset as in Section
2.4).

\propo{Proposition 9.3}{For a proper face $F$ of a polytope $P$ we
have that
\[
(\F_-(\tr_F\, P))^{op}\cong\Bl_F\, (\F_-(P))^{op}.
\]
}

Roughly speaking, a stellar subdivision makes a polytope more
simplicial, while truncation makes it more simple. This is
essential for the fact that the family of all simplicial
$d$-polytopes is dense in the family of all $d$-polytopes (see
\cite{G03}, Section 5.2, Theorem~5), and analogously, that the
family of all simple $d$-polytopes is dense in the family of all
$d$-polytopes.

\section{Stretching the interval}

Let $P$ be a simple polytope with $X$ as the set of its facets.
Then
\[
C_P=_{df}\{G^\ast\mid G\;\mbox{\rm is a nonempty face of }P\},
\]
where $G^\ast=\{F\in X\mid G\subseteq F\}$, is a simplicial
complex isomorphic to $(\F_-(P))^{op}$. (Note that the bases of
$C_P$ correspond here to the vertices of $P$.) So, for a simple
polytope $P$, our Proposition 9.3 reads
\[
C_{\tr_F P}\cong\Bl_{F\!^\ast}\, C_P.
\]
By combining this with Proposition 8.3, for $B_0$ a building set
of $C_P$, and $B_i$, for $1\leq i\leq n$ a building set of
$\Nt(C_P,B_0,\ldots,B_{i-1})$, we easily obtain that
\[
\Ss_X(C_P,B_0,\ldots, B_n)\cong C_{\tr_{F_1}(\ldots \tr_{F_m}P)},
\]
for some $F_1,\ldots, F_m$ such that $F_m$ is a proper face of
$P$, and $F_j$, $1\leq j\leq m-1$, is a proper face of
$\tr_{F_{j+1}}(\ldots \tr_{F_m}P)$. So, we have the following.

\prop{Proposition 10.1}{For $P$, $B_0,\ldots,B_n$ as above there
is a simple polytope $Q$ such that $\Ss_X(C_P,B_0,\ldots,
B_n)\cong C_Q$.}

Let $\Delta$ be an $n$-dimensional simplex with $X$ as the set of
its facets. If $B_\bot$ is the minimal building set of
$C_\Delta=P(X)-\{X\}$, i.e.\ $B_\bot=\{\{a\}\mid a\in X\}$, then
$\Ss_X(C_\Delta,B_\bot)=C_\Delta$. On the other hand, if $B_\top$
is the maximal building set of $C_\Delta$, i.e.\
$B_\top=C_\Delta-\{\emptyset\}$, then $\Ss_X(C_\Delta,B_\top)\cong
C_P$, where $P$ is an $n$-dimensional permutohedron.

By varying $B$ over all possible building sets of $C_\Delta$,
according to the definition given at the end of Section 2.4, by
relying on Proposition 3.4, we obtain the whole interval
simplex-permutohedron in dimension $n$. Proposition 9.10 of
\cite{DP10c} together with our Proposition 3.4 shows that for
every building set $B$ of $C_\Delta$ there is a simple
$n$-dimensional polytope $Q$ such that $\Ss_X(C_\Delta,B)\cong
C_Q$. According to \cite{DP10c}, $Q$ may be presented explicitly
by a finite set of inequalities (halfspaces), easily derived from
$B$. However, if one is satisfied with a less explicit
construction of $Q$, we suggest just to rely on Proposition 10.1,
which we will always do in the sequel.

It should be clear how we can modify now the interval
simplex-permuto\-he\-dron. If we replace $\Delta$ in
$\Ss_X(C_\Delta,B)$ by some other simple polytope $P$ with $X$ as
the set of its facets, and by varying $B$ over all possible
building sets of the simplicial complex $C_P$, then we obtain a
new family of simplicial complexes tied to a new interval of
simple polytopes. For example, in dimension 2, the interval
simplex-permutohedron is the interval triangle-hexagon (triangle,
quadrilateral, pentagon and hexagon). If we replace the triangle
by a quadrilateral, then we obtain the interval
quadrilateral-octagon.

For $B_\bot=\{\{a\}\mid a\in X\}$ being the minimal building set
of $C_P$, we have, as before, that $\Ss_X(C_P,B_\bot)=C_P$. So,
$P$ is the initial polytope of the new interval. On the other
hand, for $B_\top=C_P-\{\emptyset\}$ being the maximal building
set of $C_P$, by Proposition 10.1, there is a simple polytope $Q$
such that $\Ss_X(C_P,B_\top)\cong C_Q$. We may call $Q$ a
$P$-based permutohedron, and in that case the ordinary
permutohedron would be a $\Delta$-based permutohedron. So,
$P$-based permutohedron is the terminal polytope of the new
interval.

If we want to stretch the interval simplex-permutohedron and all
the other intervals obtained by replacing $\Delta$ by some other
simple polytope $P$, then it is sufficient to allow the iterated
application of $\Ss_X$. As an example, we may consider the
following four polytopes, which correspond to the simplicial
complexes described in the table of Section~7 ($X=\{x,y,z,u\}$,
$D_0=\{x+y,x+y+z\}$, $D_1=\{2x+y,2x+y+z,2x+2y+z,3x+2y+z\}$ and
$D_2=\{6x+5y+3z,3x+3y+2z\}$).

\begin{center}
\begin{picture}(300,130)(0,-10)
{\thicklines \put(50,0){\line(-1,1){50}}
\put(50,0){\line(1,1){50}} \put(50,100){\line(-1,-1){50}}
\put(50,100){\line(1,-1){50}}
\multiput(1,50)(3,0){34}{\makebox(0,0){\circle*{.5}}}
\put(50,0){\line(0,1){100}}}

\put(0,65){\vector(1,0){30}} \put(-5,65){\makebox(0,0)[r]{$x$}}

\put(100,35){\vector(-1,0){30}}
\put(105,35){\makebox(0,0)[l]{$y$}}

\put(35,0){\line(0,1){15}} \multiput(35,15)(0,3){5}{\line(0,1){1}}
\put(35,30){\vector(0,1){2}} \put(35,-5){\makebox(0,0)[t]{$z$}}

\put(65,100){\line(0,-1){15}}
\multiput(65,70)(0,3){5}{\line(0,1){1}}
\put(65,70){\vector(0,-1){2}} \put(65,105){\makebox(0,0)[b]{$u$}}

\put(-30,115){\makebox(0,0)[l]{$C=P(X)-\{X\}$}}

{\thicklines \put(230,20){\line(1,0){40}}
\put(240,30){\line(1,0){20}} \put(240,90){\line(1,0){20}}
\put(230,20){\line(-1,1){30}} \put(270,20){\line(1,1){30}}
\put(200,50){\line(1,1){40}} \put(300,50){\line(-1,1){40}}
\put(240,30){\line(0,1){60}} \put(260,30){\line(0,1){60}}
\put(230,20){\line(1,1){10}} \put(270,20){\line(-1,1){10}}}

\multiput(201,50)(3,0){34}{\makebox(0,0){\circle*{.5}}}

\put(250,100){\vector(0,-1){20}}
\put(250,105){\makebox(0,0)[b]{$x+y$}}

\put(250,0){\vector(0,1){25}}
\put(250,-5){\makebox(0,0)[t]{$x+y+z$}}

\put(180,115){\makebox(0,0)[l]{$\Et_X(C,D_0)$}}
\end{picture}
\end{center}

\begin{center}
\begin{picture}(220,175)
{\thicklines \put(50,10){\line(1,0){130}}
\put(90,40){\line(1,0){50}} \put(90,70){\line(1,0){30}}
\put(60,80){\line(1,0){20}} \put(80,150){\line(1,0){40}}

\put(50,10){\line(-1,1){20}} \put(30,30){\line(-1,2){30}}
\put(0,90){\line(6,5){60}} \put(60,140){\line(2,1){20}}
\put(120,150){\line(5,-3){100}} \put(220,90){\line(-1,-2){40}}

\put(50,10){\line(1,1){20}} \put(30,30){\line(1,1){20}}
\put(180,10){\line(-4,3){40}} \put(70,30){\line(-1,1){20}}
\put(50,50){\line(1,3){10}} \put(80,80){\line(1,-1){10}}
\put(90,70){\line(0,-1){30}} \put(70,30){\line(2,1){20}}
\put(60,80){\line(0,1){60}} \put(80,80){\line(0,1){70}}
\put(120,70){\line(0,1){80}} \put(120,70){\line(2,-3){20}}}

\multiput(2.5,90)(3,0){73}{\makebox(0,0){\circle*{.5}}}

\put(110,30){\vector(0,1){25}}
\put(110,25){\makebox(0,0)[t]{$2x+2y+z$}}

\put(70,150){\vector(0,-1){25}}
\put(70,155){\makebox(0,0)[b]{$2x+y$}}

\put(0,65){\vector(1,0){70}}
\put(-5,65){\makebox(0,0)[r]{$3x+2y+z$}}

\put(30,10){\vector(1,1){20}}
\put(27,7){\makebox(0,0)[tr]{$2x+y+z$}}

\put(-60,150){\makebox(0,0)[l]{$\Et_X(C,D_0, D_1)$}}
\end{picture}
\end{center}

\begin{center}
\begin{picture}(220,160)
{\thicklines \put(50,10){\line(1,0){130}}
\put(100,40){\line(1,0){40}} \put(100,50){\line(1,0){30}}
\put(90,70){\line(1,0){30}} \put(60,80){\line(1,0){20}}
\put(80,150){\line(1,0){40}}

\put(50,10){\line(-1,1){20}} \put(30,30){\line(-1,2){30}}
\put(0,90){\line(6,5){60}} \put(60,140){\line(2,1){20}}
\put(120,150){\line(5,-3){100}} \put(220,90){\line(-1,-2){40}}

\put(50,10){\line(1,1){20}} \put(30,30){\line(1,1){20}}
\put(180,10){\line(-4,3){40}} \put(70,30){\line(-1,1){20}}
\put(50,50){\line(1,3){10}} \put(80,80){\line(1,-1){10}}
\put(90,70){\line(0,-1){15}} \put(70,30){\line(2,1){10}}
\put(80,35){\line(1,2){10}} \put(60,80){\line(0,1){60}}
\put(80,80){\line(0,1){70}} \put(100,40){\line(0,1){10}}
\put(80,35){\line(4,1){20}} \put(90,55){\line(2,-1){10}}
\put(120,70){\line(0,1){80}} \put(120,70){\line(1,-2){10}}
\put(130,50){\line(1,-1){10}}}

\multiput(2.5,90)(3,0){73}{\makebox(0,0){\circle*{.5}}}

\put(90,30){\vector(0,1){15}} \put(93,25){\makebox(0,0)[t]{$6x\pl
5y\pl 3z$}}

\put(140,45){\vector(-1,0){25}}
\put(145,45){\makebox(0,0)[l]{$3x\pl 3y\pl 2z$}}

\put(-60,150){\makebox(0,0)[l]{$\Et_X(C,D_0, D_1, D_2)$}}

\end{picture}
\end{center}

In dimension 2, if we start with a triangle, this procedure
delivers every polygon, i.e.\ every simple two-dimensional
polytope. However, in dimension 3, if we start with a tetrahedron
it is not the case that this procedure delivers every simple
three-dimensional polytope. For example, dodecahedron would never
correspond to $\Ss_X(C_P,B)$ unless $P$ is itself a dodecahedron
and $B=B_\bot$. This is because every truncation of a simple
polytope in dimension 3 leaves at least one facet to be a triangle
or a quadrilateral.

We have described above how to produce a $P$-based permutohedron
for an arbitrary simple polytope $P$. We conclude this paper with
an example of a family of polytopes, which we call
permutohedron-based associahedra. The polytope $PA_n$ corresponds
to $\Ss_X(C_\Delta,B_\top,B)$, where $|X|=n+1$,
$C_\Delta=P(X)-\{X\}$, $B_\top=C_\Delta-\{\emptyset\}$ and $B$ is
\begin{tabbing}
\hspace{1.5em}$\{\{\{a_1,\ldots,a_k\},
\{a_1,\ldots,a_k,a_{k+1}\},\ldots,
\{a_1,\ldots,a_k,a_{k+1},\ldots,a_l\}\}\mid$
\\*
\` $a_1,\ldots,a_l$ are different elements of $X$, $1\leq k\leq
l\leq n\}$.
\end{tabbing}

Note that the bases of $\Ss_X(C_\Delta,B_\top)$ are of the form
$\alpha=\{a_1, a_1+a_2,\ldots,a_1+a_2+\ldots+a_n\}$, for
$a_1,\ldots,a_n$ different elements of $X$ and $n=|X|-1$. If we
denote $a_1+\ldots + a_k$ by $b_k$ for $1\leq k\leq n$, then
$\alpha=\{b_1,\ldots,b_n\}$ and $(P^2\sigma_X(B))_\alpha$ is of
the form $\{\{b_k, b_{k+1},\ldots, b_l\}\mid 1\leq k\leq l\leq
n\}$.

Sometimes it is much easier to present a building set by a graph
whose saturated closure (see \cite{DP10c}, Section~4) is this
building set. For example,
\begin{center}
\begin{picture}(40,20)
\put(0,5){\line(1,0){40}}

\put(0,10){\makebox(0,0)[b]{$x$}}
\put(20,10){\makebox(0,0)[b]{$y$}}
\put(40,10){\makebox(0,0)[b]{$z$}}

\put(0,5){\circle*{1.5}} \put(20,5){\circle*{1.5}}
\put(40,5){\circle*{1.5}}
\end{picture}
\end{center}
is the graph whose saturated closure is the following building set
of $P(\{x,y,z\})$
\[
\{\{x\},\{y\},\{z\},\{x,y\},\{y,z\},\{x,y,z\}\}.
\]
(This building set of $P(\{x,y,z\})$ gives rise to a 2-dimensional
associahedron, i.e.\ pentagon $K_4$.) In this sense,
$(P^2\sigma_X(B))_\alpha$ may be presented by the following graph:
\begin{center}
\begin{picture}(80,25)
\put(0,10){\line(1,0){80}}

\put(0,15){\makebox(0,0)[b]{$b_1$}}
\put(20,15){\makebox(0,0)[b]{$b_2$}}
\put(80,15){\makebox(0,0)[b]{$b_n$}}
\put(50,15){\makebox(0,0)[b]{$\ldots$}}

\put(0,10){\circle*{1.5}} \put(20,10){\circle*{1.5}}
\put(80,10){\circle*{1.5}}
\end{picture}
\end{center}
and, hence, $\Sb_X(P(\alpha),(P^2\sigma_X(B))_\alpha)$ corresponds
to an $(n-1)$-dimensional associahedron $K_{n+1}$. This means that
every vertex of the permutohedron that corresponds to
$\Ss_X(C_\Delta,B_\top)$ expands into $K_{n+1}$ in the polytope
that corresponds to
$\Ss_X(C_\Delta,B_\top,B)=\Sb(\Ss_X(C_\Delta,B_\top),P^2\sigma_X(B))$.

The bases of $\Ss_X(C_\Delta,B_\top)$ are in one to one
correspondence with the permutations of $X$ in such a way that
$a_1\ldots a_n a_{n+1}$, for $\{a_{n+1}\}=X-\{a_1,\ldots,a_n\}$,
is the permutation that corresponds to the above $\alpha$. The
bases of $\Ss_X(C_\Delta,B_\top,B)$ derived from $\alpha$ may be
interpreted as terms obtained from $a_1\ldots a_n a_{n+1}$ by
putting $n-2$ pairs of brackets (the outermost brackets are
omitted). This is done in such a way that, for example, if $n=4$,
the basis
\[
\{b_1,b_4,b_3+b_4,b_1+b_2+b_3+b_4\}
\]
of $\Ss_X(C_\Delta,B_\top,B)$, which is derived from $\alpha$, is
interpreted as $(a_1a_2)(a_3(a_4a_5))$.

Hence, $PA_n$ has the vertices, as an $n$-dimensional
permuto-associahedron $K\Pi_n$ (see \cite{Z95}, Section~9, Example
9.14) labelled by the terms built out of $n+1$ different letters
with the help of one binary operation. Some edges of $PA_n$
correspond to associativity, i.e.\ they connect two terms such
that one is obtained from the other by replacing a subterm of the
form $A\cdot(B\cdot C)$ by $(A\cdot B)\cdot C$. The other edges
correspond to transpositions of neighbours. These transpositions
are such that in a term of the form $A\cdot B$, one may permute
the rightmost letter in $A$ with the leftmost letter in $B$. So,
they correspond to ``most unexpected transposition of
neighbours''. Eventually, our family of permutohedron-based
associahedra may be taken as an alternative presentation of the
symmetric monoidal category freely generated by a set of objects.
Here is a picture of the 3-dimensional permutohedron-based
associahedron.

\begin{center}
\includegraphics[scale=0.7]{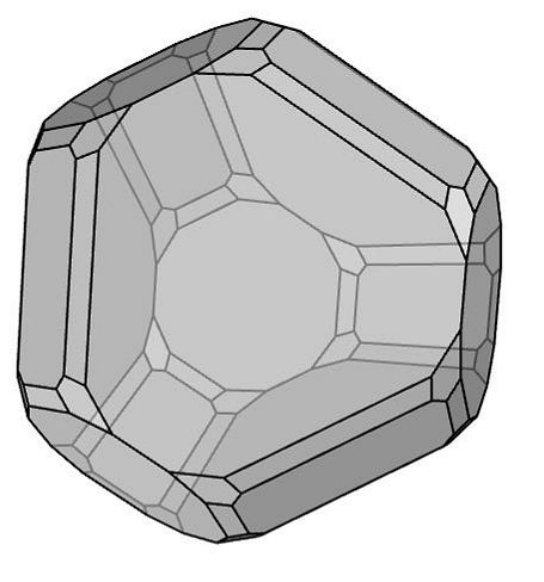}
\end{center}

\vspace{1em}

\noindent {\small \emph{Acknowledgement.} I am grateful to
Professor Jim Stasheff for some useful suggestions concerning a
previous version of this paper. I would like to thank Sonja \v
Cuki\' c for helping me to discover various possibilities of the
program package polymake (see \cite{GJ00}). Several years after
publishing this paper, Jelena Ivanovi\' c pointed out to me a
mistake in the printed picture of the 3-dimensional
permutohedron-based associahedron. I am grateful to her for this,
as well as for making the above picture of this polytope. I am
also grateful to the anonymous referees for some very useful
comments concerning the historical remarks and for finding several
misleading typos in a previous version of this paper. This work
was supported by the Ministry of Science of Serbia (Grant
ON174026). }


\begin{thebibliography}{99}

\bibitem{B89} C.\ Berge, Hypergraphs: Combinatorics of Finite Sets,
North-Holland, Amsterdam, 1989.

\bibitem{CD06} M.\ Carr,
S.L. Devadoss, Coxeter complexes and graph-associahedra, Topology
Appl. 153 (2006) 2155-2168.

\bibitem{C07} S.Lj.\ \v Cuki\' c, E.\ Delucchi,
Simplicial shellable spheres via combinatorial blowups, Proc.
Amer. Math. Soc.  135  (2007) 2403-2414.

\bibitem{DeCP95} C. De Concini, C. Procesi, Wonderful models of subspace
arrangements, Selecta Math. (N.S.) 1 (1995) 459-494.

\bibitem{D09} S.L.\ Devadoss,
A realization of graph-associahedra., Discrete Math. 309 (2009)
271-276.

\bibitem{DF08} S.L.\ Devadoss,  S.\ Forcey,
Marked tubes and the graph multiplihedron,  Algebr. Geom. Topol. 8
(2008) 2084-2108.

\bibitem{DP06} K.\ Do\v sen, Z.\ Petri\' c, Associativity as commutativity,
J. Symbolic Logic 71 (2006) 217-226 (available at: arXiv).

\bibitem{DP10a} K.\ Do\v sen, Z.\ Petri\' c,
Shuffles and concatenations in constructing of graphs, preprint,
2010 (available at: arXiv).

\bibitem{DP10b} K.\ Do\v sen, Z.\ Petri\' c,
Weak Cat-operads, preprint, 2010 (available at: arXiv).

\bibitem{DP10c} K.\ Do\v sen, Z.\ Petri\' c,
Hypergraph polytopes, Topology Appl. 158 (2011) 1405-1444
(available at: arXiv).

\bibitem{FK04} E.M.\ Feichtner, D.N.\ Kozlov, Incidence combinatorics of
resolutions, Selecta Math. (N.S.) 10 (2004) 37-60.

\bibitem{FM03} E.M.\ Feichtner, I. M\" uller, On the topology of nested set complexes, Proc.
Amer. Math. Soc. 133 (2005) 999-1006.

\bibitem{FS05} E.M.\ Feichtner, B.\ Sturmfels, Matroid polytopes, nested sets and
Bergman fans, Port. Math. (N.S.) 62 (2005) 437-468.

\bibitem{FY03} E.M.\ Feichtner, S.\ Yuzvinsky, Chow rings of toric varieties defined
by atomic lattices,  Invent. Math. 155 (2004) 515-536.

\bibitem{FS09} S.\ Forcey, D.\ Springfield,
Geometric combinatorial algebras: Cyclohedron and simplex, J.
Algebraic Combin. 32 (2010) 597-627.

\bibitem{FMP94} W.\ Fulton, R.\ MacPherson, A compactification of
configuration spaces, Ann. Math. (2) 139 (1994) 183-225.

\bibitem{Ga03} G.\ Gaiffi, Models for real subspace arrangements and
stratified manifolds, Int. Math. Res. Not. 12 (2003) 627-656.

\bibitem{Ga04} G.\ Gaiffi, Real structures of models of
arrangements, Int. Math. Res. Not. 64 (2004) 3439-3467.

\bibitem {GJ00} E.\ Gawrilow, M.\ Joswig, Polymake: a framework for
analyzing convex polytopes. Polytopes—combinatorics and
computation (Oberwolfach, 1997),  DMV Sem., 29, Birkhäuser, Basel,
2000, 43-73.

\bibitem{G03} B.\ Gr\" unbaum, Convex Polytopes, second edition, Springer, New York,
2003.

\bibitem{HRZ04} M.\ Henk, J.\ Richter-Gebert, G.M.\ Ziegler, Basic properties of convex polytopes, in J.E.\ Goodman, J.\
O'Rourke (Eds), Handbook of Discrete and Computational Geometry,
second ed., Chapman \& Hall/CRC, 2004, Section 16.

\bibitem{J07} D.\ Joji\' c, Extendable shelling, simplicial and toric h-vector of some polytopes, Publ. Inst. Math. (N.S.) 81(95) (2007) 85-93.

\bibitem{L97} J.-L.\ Loday et al. (Eds),
Operads: Proceedings of Renaissance Conferences, Contemp. Math.
202, American Mathematical Society, Providence, 1997.

\bibitem{ML98} S.\ Mac Lane, Categories
for the Working Mathematician, expanded second edition, Springer,
Berlin, 1998.

\bibitem{P09} A.\ Postnikov, Permutohedra, associahedra, and
beyond, Int. Math. Res. Not. 2009 (2009) 1026-1106.

\bibitem{PRW08} A.\ Postnikov, V.\ Reiner, L.\ Williams, Faces of
generalized permutohedra, Doc. Math. 13 (2008) 207-273.

\bibitem{S97b} J.D.\ Stasheff, From operads to physically
inspired theories (Appendix B co-authored with S.\ Shnider), in
\cite{L97}, pp.\ 53-81.

\bibitem{Z95} G.M.\ Ziegler, Lectures on Polytopes, Springer, Berlin,
1995.

\end{thebibliography}
\end{document}